\documentclass{article}

\usepackage{amssymb}
\usepackage{amsmath}
\usepackage{amsthm}

\newtheorem{theorem}{Theorem} \newtheorem{lemma}{Lemma}[section]\newtheorem{sublemma}{Sublemma}[section]
\newtheorem{propo}{Proposition}[section]

\theoremstyle{definition}
 \newtheorem{defin}{Definition}[section]

\newcommand{\ga}{\gamma}  \newcommand{\ep}{\varepsilon}
  
 \newcommand{\e}{\ep} 
 \newcommand{\vs}{\varsigma}
 \newcommand{\Z}{\mathbb{Z}}
 \newcommand{\C}{\mathbb{C}}
\newcommand{\F}{\mathbb{F}} 
\newcommand{\G} {\Gamma}

\newcommand{\Tr}{\mbox{Tr}}

\newcommand{\cN}{\mathcal{N}}
\newcommand{\cH}{\mathcal{H}}
\newcommand{\cA}{\mathcal{A}} 
\newcommand{\Sub}{\mbox{Sub}}\newcommand{\Orb}{\mbox{Orb}}
\newcommand{\Stab}{\mbox{Stab}}\newcommand{\Sch}{\mbox{Sch}}
\newcommand{\vi}{\vskip 0.1in \noindent}
\newcommand{\roo}{\mbox{root}}
\newcommand{\Forg}{\mbox{F}}
\newcommand{\dist}{\mbox{dist}}
\newcommand{\diam}{\mbox{diam}}\newcommand{\Cay}{\mbox{Cay}}
\newcommand{\sk}{\Sch^{K,Q}_\G}
\title{Uniformly recurrent subgroups and simple $C^*$-algebras\footnote{
by the ERC Consolidator Grant "Asymptotic invariants of discrete groups,
sparse graphs and 
locally symmetric spaces" No. 648017. } }

\author{G\'abor Elek}

\begin{document}
\maketitle

\begin{abstract}
We study uniformly recurrent subgroups (URS) introduced by Glasner and Weiss
\cite{GW}. Answering their query we show that any URS $Z$
of a finitely generated group is the stability system of a minimal $Z$-proper
action. We also show that for any sofic $URS$ $Z$ there is a $Z$-proper action
admitting an invariant measure. We prove that for a $URS$ $Z$ all
$Z$-proper actions admits an invariant measure if and only if $Z$ is
coamenable. In the second part of the paper we study the separable $C^*$-algebras
associated to URS's. We prove that if a URS is generic then its
$C^*$-algebra is simple. We give various examples of generic URS's with exact and
nuclear $C^*$-algebras and an example of a URS $Z$ for which the
associated
simple $C^*$-algebra is not exact and not even locally reflexive, in particular, it admits both a
uniformly amenable trace and a nonuniformly amenable trace. 
\end{abstract}

\newpage
\tableofcontents
\newpage
\section{Introduction}
Let $\G$ be a countable group and $\Sub(\G)$ be the compact space of
all subgroups of $\G$. The group $\G$ acts on $\Sub(\G)$ by conjugation.
{\bf Uniformly recurrent subgroups}\,(URS) were defined by 
Glasner and Weiss \cite{GW}
as closed, invariant subsets $Z\subset \Sub(\G)$ such that the action
of $\G$ on $Z$ is minimal (every orbit is dense).
Now let $(X,\G,\alpha)$ be a $\G$-system (that is, $X$ is a compact metric
space and $\alpha:\G\to \mbox{Homeo}(X)$ is a homomorphism).
For each point $x\in X$ one can define the topological stabilizer subgroup 
$\Stab_\alpha^0(x)$ by
$$\Stab_\alpha^0(x)=\{\ga\in \G\,\mid\, \ga \,\mbox {fixes some neighborhood of $x$}\}\,.$$

\noindent
Let us consider the  $\G$-invariant subset $X^0\subseteq X$ such that
$x\in X^0$ if and only if $\Stab_\alpha(x)=\Stab_\alpha^0(x)\,.$  The closure of the
invariant subset 
$\Stab_{\alpha} (X^0)\subset \Sub(\G)$ is
called the {\bf stability system} of $(X,\G,\alpha)$ 
(see also \cite{Kaw},\cite{BB}).
If the action is minimal, then the stability system of $(X,\G,\alpha)$
is a URS.
Glasner and Weiss proved (Proposition 6.1,\cite{GW}) that
for every URS $Z\subset \Sub(\G)$ there exists a topologically transitive
(that is there is a dense orbit) system $(X,\G,\alpha)$ with $Z$ as its
stability system.
They asked (Problem 6.2., \cite{GW}), whether for any URS $Z$
there exists a 
minimal system $(X,\G,\alpha)$ with $Z$ as its stability system. Recently,
Kawabe \cite{Kaw} gave an affirmative answer for this question in the case
of amenable groups. 
\begin{defin}\label{proper}
Let $\G$ be a countable group
and $Z\subset \Sub(\G)$ be a URS. A system  $(X,\G,\alpha)$ is $Z$-proper
if for any $x\in X$ $\Stab_\alpha(x)=\Stab_\alpha^0(x)$ and $\Stab_\alpha(X)\in Z$.
\end{defin}
\noindent
Before stating our first result we prove a lemma for the sake of completeness.
\begin{lemma} \label{simplelemma}
If $(X,\Gamma,\alpha)$ is a $Z$-proper system, then the map $\Stab_\alpha:X\to \Sub(\Gamma)$ is
continuous. 
\end{lemma}
\proof
Let $\{x_n\}^\infty_{n=1}$ be a sequence in $X$ converging to an element $x\in X$. We need to
show that $\gamma\in\Stab_\alpha(x)$ if and only if there exists some constant $N_\gamma>0$ such that
if $n\geq N_\gamma$ then $\gamma\in\Stab_\alpha(x_n)$. Clearly, if $\gamma\in\Stab_\alpha(x_n)$, then by
the continuity of $\alpha$, $\gamma\in\Stab_\alpha(x)$ for large enough $n$. In other words, for any
$\G$-system $(X,\Gamma,\alpha)$ the map $\Stab_\alpha:X\to \Sub(\Gamma)$ is upper-semicontinuous. It is
important to note that for $\Gamma$-systems in general the map $\Stab_\alpha$ is not necessarily continuous at all
points $x\in X$. Let $(X,\alpha,\Z)$ be the standard Bernoulli shift. That is, $X=\{0,1\}^\Z$ and $\alpha$ is the 
left translation by $\Z$. Let $x_n\in X$ be defined the following way. For $n\geq 1$, let $x_n(k)=1$ if $|k|\leq n$ and let
$x_n(k)=0$, otherwise. Also, let $x(k)=1$ for any $k\in\Z$. Then, $x_n\to x$. On the other hand,
$\Stab_\alpha(x_n)=\{0\}$ for all $n\geq 1$ and $\Stab_\alpha(x)=\Z$. Now, if $(X,\alpha, \Gamma)$ is $Z$-proper for some
URS $Z$ and $\gamma\in\Stab_\alpha(x)$, then $\gamma\in\Stab_\alpha(y)$ for some neighborhood $x\in U \subset X$. Hence,
there exists some constant $N_\gamma>0$ such that $\gamma\in Stab_\alpha(x_n)$ provided that $n\geq N_\gamma$. Therefore our
lemma follows. \qed
\begin{theorem}\label{tetel1}
If $\G$ is a finitely generated group and $Z\subset \Sub(\G)$ is a URS, then
there exists a minimal $Z$-proper system $(X,\G,\alpha)$ (that is,  $Z$ is the stability
system of  $(X,\G,\alpha)$).
\end{theorem}
\noindent
In fact, we will show that $X$ can be chosen as a $Z$-proper minimal Bernoulli
subshift (see Definition  \ref{zbern}).
In the proof we will use the Lov\'asz Local Lemma technique of 
Alon, Grytczuk, Haluszczak and Riordan
 \cite{AGHR} to construct a minimal action
on the space of rooted colored $\G$-Schreier graphs. 
This approach has already
been used to construct free $\G$-Bernoulli subshifts by Aubrun, Barbieri and Thomass\'e 
\cite{ABT} . Very recently, Matte Bon and Tsankov \cite{MT} completely answered the
query of Glasner and Weiss for uniformly recurrent subgroups of discrete and
locally compact groups.
The next result of the paper is about the existence of invariant measures on
$Z$-proper Bernoulli subshifts. For a long time all finitely generated groups that
had been known to have free Bernoulli subshifts were residually-finite.
Then Dranishnikov and Schroeder \cite{DS} constructed a free Bernoulli subshift
for any torsion-free hyperbolic group. Somewhat later Gao, Jackson and Seward
proved that any countable group has free Bernoulli subshifts \cite{GJS1}, 
\cite{GJS2}. 
On the other hand, Hjorth and Molberg \cite{HM} proved that for any countable group $\Gamma$ there exists a free continuous action of $\Gamma$
on a Cantor set admitting an invariant probability measure. We will prove the following
result.
\begin{theorem} \label{tetel2}
Let  $\G$ be a finitely generated group
and $Z\subset\Sub(\Gamma)$ be a sofic URS (see Definition \ref{sofurs}) then
there exists a $Z$-proper Bernoulli shift with an invariant probability measure. In
particular, for every finitely generated sofic group $\G$ there exists a free
Bernoulli subshift with an invariant probability measure.
\end{theorem}
\noindent
 Immediately after the first
version of our paper appeared, using a measurable version of the Local Lemma,
Bernhsteyn \cite{Bern} proved that free Bernoulli
subshift admitting an invariant probability measure exists for any countable group. He also noted that this result
follows from a deep theorem of Seward and Tucker-Drob \cite{ST}.
We can actually characterize those uniformly recurrent subgroups $Z$ for which
all the $Z$-proper actions admit invariant probability measures (Theorem
\ref{karak}).
\vskip 0.2in
\noindent
The second part of the paper is about $C^*$-algebras. For any finitely generated group  $\G$ and uniformly recurrent subgroup $Z\subset\Sub(G)$, we associate
a separable $C^*$-algebra $C^*_r(Z)$. For any group $\G$, if $Z=\{1\}$, the associated
$C^*$-algebra $C^*_r(Z)$ is just the reduced $C^*$-algebra of the group. It is
known that the reduced $C^*$-algebra of a group $\G$ is simple if and only if
the group admits no non-trivial amenable uniformly recurrent subgroups
\cite{Ken}.
We prove (Theorem \ref{ursimple}) that if the URS $Z$ is generic (see Subsection \ref{genericity})  then the
$C^*$-algebra $C^*_r(Z)$  is always simple. Using the coloring scheme developed in the first part of the paper,  we will show how to
construct generic URS's from a single infinite graph of bounded vertex
degrees. By this construction we obtain examples of generic URS's with nuclear
(Theorem \ref{manuclear})
and exact( but not nuclear) $C^*$-algebras (Proposition \ref{maexact}).
Finally, we will construct a generic $URS$
$Z$ for which the simple $C^*$-algebra $C^*_r(Z)$ is not locally-reflexive
(hence not exact). In fact, this algebra $C^*_r(Z)$ admits both a uniformly amenable and a non-uniformly
amenable trace. We will see that the URS above is not Borel equivalent to a
free minimal action of any countable group. 
\section{Schreier graphs}
\subsection{The space of rooted Schreier graphs}\label{rooted}
Let $\Gamma$ be a finitely generated group with a 
generating system $Q=\{\gamma_i\}^n_{i=1}$. 
Let $H\in\Sub(\Gamma).$ Then the Schreier graph associated to  $H$  is 
constructed as follows.
\begin{itemize}
\item The vertex set of the Schreier graph of $H$ is the coset space $\Gamma/H$ (that is the group
$\Gamma$ acts on the vertex set of the Schreier graphs on the left).
\item The vertices corresponding to the cosets $aH$ and $bH$ are connected
by a directed edge labeled by the generator $\gamma_i$ if $\gamma_iaH=bH$, or by $\gamma_i^{-1}$ if
$\gamma_i bH=aH\,$ (note that we allow loops and multiply labeled directed edges).
\end{itemize}
\vi
The coset class of $H$ is called the {\it root} of the Schreier graph associated to $H$. The set
of all rooted Schreier graphs will be denoted by $\Sch^Q_\Gamma$. So, we have a map $S^Q_\Gamma:\Sub(\G)\to \Sch^Q_\Gamma$
such that $S^Q_\G(H)$ is the rooted Schreier graph associated to the subgroup $H$.
We will consider the usual shortest path distance on the graph $S^Q_\G(H)$ and
denote the ball of radius $r$ around the root $H$ by $B_r(S^Q_\G(H), H)\,.$
Note that  $B_r(S^Q_\G(H), H)$ is a rooted edge-labeled graph.
The space of all Schreier graphs $\Sch^Q_\Gamma$ is a compact metric
space, where
$$d_{\Sch^Q_\G}(S^Q_\G(H_1), S^Q_\G(H_2))=2^{-r}\,,$$
if $r$ is the largest integer for which the $r$-balls
$B_r(S^Q_\G(H_1), H_1)$ and \\ $B_r(S^Q_\G(H_2), H_2)$ are rooted-labeled isomorphic. 
We can define the action of the group $\Gamma$ on the compact metric space $\Sch^Q_\Gamma$ in the
following way.
If $\ga\in\Gamma$ and $H\in\Sub(\Gamma)$,
then 
$$\ga(S^Q_\Gamma(H))=S^Q_\Gamma(\ga H \ga^{-1})\,.$$
The graph 
$S^Q_\Gamma(\ga H \ga^{-1})$ can be regarded as the same graph
as $S^Q_\Gamma(H)$ with the new root $\ga H$.
We will use the root-change picture of the $\Gamma$-action on
$\Sch^Q_\G$ later in the paper. If $S=S^Q_\Gamma(H)$ is a Schreier graph  and $x=\ga H$ is another vertex of $S$, then
$(S^Q_\G(H),x)$ will denote the Schreier graph with underlying labeled graph $S$ and root $x$. In this case $(S^Q_\G,x)$ is isomorphic
to $S^Q_\G(\ga H\ga^{-1})$ as rooted Schreier graphs.
Clearly, $S^Q_\G:\Sub(\Gamma)\to \Sch^Q_\Gamma$, is a 
homeomorphism commuting
with the $\G$-actions defined above. Let $Z\subset \Sub(\Gamma)$, then $S^Q_\G(Z)\subset \Sch^Q_\G$ is
a closed $\Gamma$-invariant subspace of rooted Schreier graphs.  
\subsection{Schreier graphs and uniformly recurrent subgroups}
\begin{propo}\label{repet}
Let $\G$ and $Z$ be as above, $H\in Z$ and $S^Q_\G(H)$ be the corresponding rooted Schreier graph. Then for any $x\in V(S^Q_\G(H))$ and $R>0$ 
there exists
$S_{x,R}>0$ such that for any $y\in V(S^Q_\G(H))$, there is a $z\in V(S^Q_\G(H))$ so that
\begin{itemize}
\item $d_{S^Q_\G(H)}(y,z)\leq S_{x,R}$
\item The rooted labeled balls $B_R(S^Q_\G(H),x)$ and $B_R(S^Q_\G(H),z)$ are isomorphic.
\end{itemize}

\noindent
Conversely, if $H\in\Sub(\G)$ has the repetition property as above, then its orbit closure in $\Sub(\G)$ is a uniformly
recurrent subgroup.
\end{propo}
\proof
We proceed by contradiction. Suppose that there is some $x\in V(S^Q_\G(H))$ such that for all $n\geq 1$ there exists
$y_n\in V(S^Q_\G(H))$ such that if $d_{S^Q_\G(H)}(y,z)\leq n$, then $B_R(S^Q_\G(H),x)$ and $B_R(S^Q_\G(H),z)$ are not isomorphic.
Let $S\in \Sch^Q_\G$ be a rooted Schreier graph that is a limitpoint of the sequence of rooted Schreier graphs $\{S^Q_\G(H),y_n\}^\infty_{n=1}$.
Then, if $q\in V(S)$, the rooted balls  $B_R(S^Q_\G(H),x)$ and $B_R(S,q)$ are not isomorphic. Hence, the orbit closure of $S$ in the $\Gamma$-space
$\Sch^Q_\G$ does not contain the Schreier graph $S^Q_\G(H)$ in contradiction with the minimality of $Z$.

\noindent
Now we prove the converse. Let $H\in \Sub(\Gamma)$ be a subgroup satisfying the condition
of our lemma. Let $K,L\in\Sub(\Gamma)$ be elements
of the orbit closure of $H$. It is enough to show that the orbit closure of $K$ contains $L$.
Let $R>0$ be an integer. We need to show that there exists
$x\in V(S^Q_\G(K))$ such that
$B_R(S^Q_\G(K),x)$ is rooted-labeled isomorphic to $B_R(S^Q_\G(L),L))\,.$
Since $L$ is in the orbit closure of $H$, we have $y\in V(S^Q_\G(H))$ such that
$B_R(S^Q_\G(H),y)$ is rooted-labeled isomorphic to $B_R(S^Q_\G(L),L))\,.$
By our condition, if $K$ is in the orbit closure of $H$, there exists
$x\in V(S^Q_\G(K))$ so that $B_R(S^Q_\G(K),x)$ is rooted-labeled isomorphic to
$B_R(S^Q_\G(H),y)$. This finishes the proof of our proposition. \qed
\subsection{Genericity}\label{genericity}
Let $\G$ be as above and $Z\subset \Sub(\Gamma)$ be a URS.
We say that $Z$ is {\bf generic} if for every $H\in Z$, the coset space $\G/H$ and the
orbit of $H$ in $\Sub(\G)$ are $\G$-isomorphic sets under the map $\phi:\Gamma/H\to \mbox{Orb}(H)$, $\phi(gH)=gHg^{-1}$. That is, all the elements of $Z$ are self-normalizing subgroups. We will
give several examples of generic URS's in Section \ref{example}.
\begin{propo}
Let $Z$ be a generic URS of $\G$. Then for each $H\in Z$, $\Stab^0_\alpha(H)=\Stab_\alpha(H)=H\,.$ That is,
$(Z,\Gamma,\alpha)$ is a $Z$-proper system, where $\alpha$ is the conjugation action of
$\G$ on $Z$. Hence, the stability
system of a generic $URS$ is itself.
\end{propo}
\proof Let $H\in Z$. Then by genericity, $\Stab_\alpha(H)$ is the stabilizer of the root in $S^Q_\G(H)$, that is,
$\Stab_\alpha(H)=H.$
Also, if $h\in H$, then $h$ fixes the root of every element of $\Sch^Q_\G$ that is close enough to
$S^Q_\G(H)$, hence $\Stab^0_\alpha(H)=\Stab_\alpha(H)\,.$ \qed
\begin{propo}\label{assy}
The uniformly recurrent subgroup $Z$ is generic if and only if the following statement holds.
For any $R>0$ there exists $S>0$ such that if $H\in Z$, $x,y\in V(S^Q_\G(H))$, $0<d_{S^Q_{\G}(H)}(x,y)\leq R$, then
the rooted balls $B_S(S^Q_\G(H),x)$ and $B_S(S^Q_\G(H),y)$ are not rooted-labeled isomorphic.
\end{propo}
\proof
Suppose that for any $n\geq 1$, there exists $H_n\in Z$ and $x_n,y_n\in\G/H_n$
such that
\begin{itemize}
\item
$0<d_{\G/ H_n}(x_n,y_n)\leq R\,.$
\item
the $n$-balls around $x_n$ and $y_n$ are rooted-labeled isomorphic.
\end{itemize}

\noindent
Let $(S^Q_\G(H),H)$ be a limitpoint of the sequence
$\{(S^Q_\G(H_n),x_n)\}^\infty_{n=1}$ in $S^Q_\G$. Then, there exists $\gamma\in \G,\gamma\notin H$, so that
$(S^Q_\G(H),H)$ and $(S^Q_\G(H),gH)$ are rooted-labeled isomorphic. Hence
$\phi:\G/H\to\mbox{Orb}(H)$ is not a bijective map.
On the other hand, it is clear that if the condition of our proposition
is satisfied for any $H\in Z$, then $\phi:\G/H\to \mbox{Orb}(x)$ is always 
bijective, hence $Z$ is generic. \qed

\subsection {The Bernoulli shift space of uniformly recurrent subgroups}
Let $\G,Q$ be as in the previous subsection, $H\in\Sub(\Gamma)$ and let $K$ be a finite alphabet. 
A rooted $K$-colored
Schreier graph of $H$ is the rooted Schreier graph $S^Q_\G(H)$ equipped
with a vertex-coloring $c:\Gamma/H\to K$. 
Let $\sk$ be the set of all rooted $K$-colored Schreier-graphs. Again,
we have a compact, metric topology on $\sk$:
$$d_{\sk}(S,T)=2^{-r}\,,$$
\vi
if $r$ is the largest integer such that the $r$-balls around the
roots of the graphs $S$ and $T$ are rooted-colored-labeled isomorphic.
We define $d_{\sk}(S,T)=2$ if the $1$-balls around the roots are nonisomorphic
and even the colors of the roots are different.
Again, $\Gamma$ acts on the compact space $\sk$ by the root-changing map.
Hence, we have a natural color-forgetting map $\Forg:\sk\to \Sch^Q_\Gamma$  
that commutes with the $\Gamma$-actions. Notice that if
a sequence $\{S_n\}^\infty_{n=1}\subset \sk$ converges to $S\in\sk$, then 
for any $r\geq 1$ there exists some integer $N_r\geq 1$ such that
if $n\geq N_r$ then the $r$-balls around the roots of the graph $S_n$ and
the graph $S$ are rooted-colored-labeled isomorphic. 
Let $H\in\Sub(\Gamma)$ and $c:\G/H\to K$ be a vertex coloring
that defines the element  $S_{H,c}\in \sk$. Then of course,
$\gamma(S_{H,c})=S_{H,c}$ if $\gamma\in H$. On the other hand, if
$\gamma(S_{H,c})=S_{H,c}$ and $\gamma\notin H$ then we have the
following lemma that is immediately follows from the definitions of the $\G$-actions.
\begin{lemma}\label{trivi}
Let $\gamma\notin H$ and $\gamma(S_{H,c})=S_{H,c}$. Then there exists
a colored-labeled graph-automorphism of the $K$-colored
labeled graph $S_{H,c}$ moving
the vertex representing  $H$ to the vertex representing $\gamma(H)\neq H$.
\end{lemma}
\vi
Note that we have a continuous $\G$-equivariant map $\pi:\sk\to\Sub(\Gamma)$,
where $\pi(t)=(S^Q_\G)^{-1}\circ \Forg(t).$
Let $Z$ be a URS of $\Gamma$.  We say that the element
$t\in \sk$ is {\bf $Z$-regular} if $\pi(t)=H\in Z$ and  $\Stab_\alpha(t)=H$, where $\alpha$ is the
left action of $\Gamma$ on $\sk$. Note that if $H\in Z$ and $t$ is a $K$-coloring
the Schreier graph $S^Q_\G(H)$, then by Lemma \ref{trivi}, $t$ is $Z$-regular if and only
if there is no non-trivial colored-labeled automorphism of $t$.
\begin{propo}\label{propo1}
Let $Y\subset \sk$ be a closed $\G$-invariant subset consisting
of $Z$-regular elements. Let
$(M,\G,\alpha)\subset (Y,\G,\alpha)$ be a minimal $\Gamma$-subsystem.
Then $M$ is $Z$-proper, that is,  for any $m\in M$, $\Stab_\alpha^0(m)=\Stab_\alpha(m)\in Z$. Also,
$\pi(M)=Z$.
\end{propo}
\proof
Let $h\in \Stab_\alpha(m)$. Then by $Z$-regularity $h$ fixes the root of $m$.
Therefore, $h$ fixes the root of $m'$ provided that $d_{\sk}(m,m')$ is small
enough. Thus, $h\in \Stab_\alpha^0(m)$. Since $\pi$ is a
$\Gamma$-equivariant continuous map and $M$ is a closed $\Gamma$-invariant
subset, $\pi(M)=Z$. \qed

\begin{defin} \label{zbern}
Let $Z$ be as above and $K$ be a finite alphabet. Let $B^K(Z)$ be the
$\G$-invariant subset of all elements $S$ of $\Sch_\G^{K,Q}$ such that the
underlying Schreier graph is in $\Sch^Q_\G(Z)$.
We call $B^K(Z)$ the $K$-Bernoulli shift space of $Z$. A closed
$\G$-invariant subset of $B^K(Z)$ is called a Bernoulli subshift of $Z$.
\end{defin}

\vskip 0.1in
\noindent
Note that if $Z=\{1\}$, then $Z$-properness is just the classical notion
of $\G$-freeness, and the $Z$-subshifts are the Bernoulli subshifts
of $\G$.

\section{Lov\'asz's Local Lemma and the proof of Theorem \ref{tetel1}}
Let $Z$ be a URS of $\G$.
By Proposition \ref{propo1}, it is enough to construct a closed
$\Gamma$-invariant subset $Y\subset \sk$ for some alphabet $K$ such
that all the elements of $Y$ are $Z$-proper. This will give us a bit more than just a continuous action
having stability system $Z$, $Y$ will be a minimal $Z$-proper Bernoulli subshift. It is quite clear that
the stability system of a minimal $Z$-proper Bernoulli subshift is always $Z$ itself.
Let $H\in Z$ and consider the Schreier graph $S=S^Q_\G(H)$. Following \cite{ABT} and
\cite{AGHR} we call a coloring $c:\Gamma/H\to K$ nonrepetitive if
for any path $(x_1,x_2,\dots, x_{2n})$ in $S$ there exists some
$1\leq i\leq n$ such that
$c(x_i)\neq c(x_{n+i})\,.$ We call all the other colorings repetitive.
\begin{theorem}\label{alon}[Theorem 1 \cite{AGHR}]
For any $d\geq 1$ there exists a constant $C(d)>0$ such that any 
graph $G$ (finite or infinite) with vertex degree bound $d$ has
a nonrepetitive coloring with an alphabet $K$, provided that  $|K|\geq C(d)$.
\end{theorem}
\proof Since the proof in \cite{AGHR} is about edge-colorings and
the proof in \cite{ABT} is in slightly different setting, for completeness
we give a proof using  Lov\'asz's  Local Lemma, that closely follows the
proof in \cite{AGHR}. Now, let us state the Local Lemma.
\begin{theorem}[The Local Lemma]
Let $X$ be a finite set and $\Pr$ be a probability distribution on
the subsets of $X$. For $1\leq i\leq r$ let $\cA_i$ be a set of 
events, where an ``event'' is just a subset of $X$. Suppose that for
all $A\in \cA_i$, $\Pr(A_i)=p_i$.
 Let $\cA=\cup^r_{i=1} \cA_i$. Suppose that there are real numbers
$0\leq a_1,a_2,\dots,a_r<1$  and $\Delta_{ij}\geq 0$, $i,j=1,2,\dots, r$
such that the following conditions hold:
\begin{itemize}
\item for any event $A\in  \cA_i$ there exists a set $D_A\subset \cA$ with
$|D_A\cap \cA_j|\leq \Delta_{ij}$ for all $1\leq j \leq r$ such that
$A$ is independent of $\cA\backslash (D_A\cup \{A\}),$
\item $p_i\leq a_i\prod^r_{j=1}(1-a_j)^{\Delta_{ij}}$ for all $1\leq i \leq r\,.$
\end{itemize}
Then $\Pr(\cap_{A\in\cA} \overline{A})>0$. 
\end{theorem}
\vi
Let $G$ be a finite graph with maximum degree $d$. It is enough to
prove our theorem for finite graphs. Indeed, if $G'$ is a connected
infinite graph with vertex degree bound $d$, then for each ball around
a given vertex $p$ we have a nonrepetitive coloring. Picking a pointwise
convergent subsequence of the colorings we obtain a nonrepetitive 
coloring of our
infinite graph $G'$. 

\noindent
Let $C$ be a large enough number, its exact value will
be given later. Let $X$ be the set of all random $\{1,2,\dots,C\}$-colorings
of $G$. Let $r=\diam(G)$ and for $1\leq i \leq r$ and for any path $P$ of 
length $2i-1$ 
let $A(P)$ be the event that $P$ is repetitive. Set
$$\cA_i=\{A(P):\, \mbox{$P$ is a path of length $2i-1$ in $G$}\}\,.$$
\vi
Then $p_i=C^{-i}$. The number of paths of length $2j-1$ that intersects
a given path of length $2i-1$ is less or equal than $4ij d^{2j}$.
So, we can set $\Delta_{ij}=4ijd^{2j}$.
Let $a_i=\frac{1}{2^id^{2i}}$. Since $a_i\leq \frac{1}{2}$, we have that
$(1-a_i)\geq \exp(-2a_i)$.
In order to be able to apply the Local Lemma, we need
that for any $1\leq i \leq r$
$$p_i\leq a_i\prod_{j=1}^r \exp(-2a_j\Delta_{ij})\,.$$
That is
$$C^{-i}\leq a_i\prod _{j=1}^r \exp(-8ija_jd^{2j})\,,$$
or equivalently
$$C\geq 2d^2\exp\left( 8\sum_{j=1}^r \frac{j}{2^j}\right)\,.$$
Since the infinite series $\sum_{j=1}^\infty \frac{j}{2^j}$ converges to $2$,
we obtain that for large enough $C$, the conditions of the Local Lemma
are satisfied independently on the size of our finite graph $G$. This ends
the proof of Theorem \ref{alon}.\quad \qed
\vi
Let $|K|=C(|Q|)$ and let $c:\Gamma/H\to K$ be a nonrepetitive $K$-coloring
that gives rise to an element $y\in \sk$. The following
proposition finishes the proof of Theorem \ref{tetel1}.
\begin{propo} \label{propo2}
All elements of the orbit closure $Y$ of $y$ in $\sk$ are $Z$-regular.
\end{propo}
\proof Let $x\in Y$ with underlying Schreier graph $S^Q_\G(H')$ and
coloring $c':\Gamma/H' \to K$. Since $Z$ is a URS, $H'\in Z$. Indeed,
$\pi^{-1}(Z)$ is a closed $\Gamma$-invariant set and $y\in \pi^{-1}(Z)$.
Clearly, $\alpha(\gamma)(x)=x$ if $\gamma\in H'$. Now suppose that
$\alpha(\gamma)(x)=x$ and $\gamma\notin H'$ (that is $x$ is not $Z$-proper).
By Lemma \ref{trivi}, there exists a colored-labeled automorphism $\theta$
of the graph $x$ moving $\roo(x)$ to $\alpha(\gamma)(\roo(x))\neq \roo(x)$. 
Note that if $a$ is a vertex of $x$, then $\theta(a)\neq a$. Indeed, if a labeled
automorphism of a Schreier graph fixes one vertex, it must fix all the other vertices as well.
Now we proceed similarly as in the proof of Lemma 2 \cite{AGHR} or 
in the proof of Theorem 2.6 \cite{ABT}.
Let $a\in V(x)$ be a vertex such that
there is no $b\in x$ such that
$\dist_x(b,\theta(b))<\dist_x(a,\theta(a))$.
Let $(a=a_1, a_2,\dots, a_{n+1}=\theta(a))$ be a shortest path between $a$ and
$\theta(a)$.
For $1\leq i \leq n$, let
$\alpha(\gamma_{k_i})(a_i)=a_{i+1}\,.$ Then let $a_{n+2}=\alpha(\gamma_{k_1})(a_{n+1}),
a_{n+3}=\alpha(\gamma_{k_2})(a_{n+2}),\dots, a_{2n}=\alpha(\gamma_{k_n})(a_{2n-1})\,.$
Since $\theta$ is a colored-labeled automorphism, for any $1\leq i \leq n$\,
\begin{equation}\label{e1}
c(a_i)=c(a_{i+n})\,.
\end{equation}
\begin{lemma} 
The walk $(a_1,a_2,\dots, a_{2n})$ is a path.
\end{lemma}
\proof Suppose that the walk above crosses itself, that is for some $i,j$,
$a_j=a_{n+i}$. If $(n+1)-j\geq (n+i)-(n+1)=i-1\,,$ then
$\dist_x(a_2,\theta(a_2))=\dist_x(a_2,a_{n+2})<\dist_x(a,\theta(a))\,.$
On the other hand, if $(n+1)-j\leq (n+i)-(n+1)=i-1\,,$ then
$\dist_x(a_n,\theta(a_n))=\dist_x(a_n,a_{2n-1})<\dist_x(a,\theta(a))\,.$
Therefore, $(a_1,a_2,\dots, a_{2n})$ is a path. \qed
\vi
By (\ref{e1}) and the previous lemma, the $K$-colored Schreier-graph $x$
contains a repetitive path. Since $x$ is in the orbit closure of $y$, this
implies that $y$ contains a repetitive path as well, in contradiction with
our assumption. \qed

\section {Sofic groups, sofic URS's and invariant measures}
\subsection{Sofic groups}
First, let us recall the notion of a finitely generated sofic group. 
 Let $\G$ be a finitely generated
infinite group with a symmetric generating system 
$Q=\{\gamma_i\}^n_{i=1}$ and
a surjective homomorphism $\kappa:\F_n\to \Gamma$ from the
free group $\F_n$ with generating system $\overline{Q}=\{r_i\}^n_{i=1}$ mapping
$r_i$ to $\gamma_i$. Let $\Cay^Q_\G$ be the Cayley graph of $\G$ with
respect to the generating system $Q$, that is the Schreier graph corresponding
to the subgroup $H=\{1_\G\}$. Let $\{G_k\}^\infty_{k=1}$ be a sequence
of finite $\F_n$-Schreier graphs.
We call a vertex $p\in V(G_k)$ a $(\Gamma,r)$-vertex if there exists
a rooted isomorphism
$$\Psi:B_r(G_k,p)\to B_r(\Cay^Q_\G,1_\G)\,$$
\vi
such that if $e$ is a directed edge in the ball $B_r(G_k,p)$ labeled
by $r_i$, then the edge $\Psi(e)$ is labeled by $\gamma_i$.
We say that $\{G_k\}^\infty_{k=1}$ is a sofic approximation of $\Cay^Q_\G$, if
for any $r\geq 1$ and a real number $\e>0$ there exists $N_{r,\e}\geq 1$
such that if $k\geq N_{r,\e}$ then there exists a subset $V_k\subset V(G_k)$
consisting of $(\Gamma,r)$-vertices such that $|V_k|\geq (1-\e)|V(G_k)|$.
A finitely generated group $\Gamma$ is called sofic if the Cayley-graphs of $\G$ admit
sofic approximations. Sofic groups were introduced by Gromov in \cite{Gro} under the name
 of initially subamenable groups, the word ``sofic'' was coined by Weiss in
\cite{Weiss}. It is important to note that all the amenable, residually-finite
and residually amenable groups are sofic, but there exist finitely
generated sofic groups that are not residually amenable (see the book
of Capraro and Lupini \cite{Lup} on sofic groups). It is still an open
question whether all groups are sofic.
\subsection{Sofic URS's and the proof of Theorem \ref{tetel2}}\label{sub42}
We can extend the notion of soficifty from groups to URS's in the following way.
Let $\G$,$Q$, $\kappa:\F_n\to \G$ be as above and let $Z\subset \Sub(\G)$ be a uniformly recurrent
subgroup. Again, let $\{G_k\}^\infty_{k=1}$ be a sequence
of finite $\F_n$-Schreier graphs. We call a vertex $p\in V(G_k)$ be a $(Z,r)$-vertex if there exists
a rooted isomorphism
$\Psi:B_r(G_k,p)\to B_r(S^Q_\G(H),H)$, where $H\in Z$ such that
if $e$ is a directed edge in the ball $B_r(G_k,p)$ labeled
by $r_i$, then the edge $\Psi(e)$ is labeled by $\gamma_i$. Similarly to the case of groups, we say
that $\{G_k\}^\infty_{k=1}$ is a sofic approximation of the uniformly recurrent subgroup $Z$ if
or any $r\geq 1$ and a real number $\e>0$ there exists $N_{r,\e}\geq 1$
such that if $k\geq N_{r,\e}$ then there exists a subset $V_k\subset V(G_k)$
consisting of $(Z,r)$-vertices such that $|V_k|\geq (1-\e)|V(G_k)|$. 
\begin{defin}\label{sofurs} A uniformly recurrent subgroup
is  sofic if it admits a sofic approximation system (note that soficity does not depend
upon the choice of the generating system)
\end{defin}

\noindent
In Section \ref{example} we will construct a large variety of generic and non-generic URS's.
The rest of this subsection is devoted to the proof of Theorem \ref{tetel2}.
Let $Z$ be a sofic URS and $\{G_k\}^\infty_{k=1}$ be a sofic approximation of $Z$.
Using Theorem \ref{alon}, for each $k\geq 1$ let us choose a nonrepetitive coloring
$c_k:V(G_k)\to K$, where $|K|\geq C(|Q|)$. We can associate
a probability measure $\mu_k$ on the space
of $K$-colored $\F_n$-Schreier graphs $\Sch^{\overline{Q},K}_{\F_n}$, where
$\overline{Q}=\{r_i\}^n_{i=1}$ is the generating system of the free group $\F_n$. Note that
the origin of this construction can be traced back to the paper of Benjamini
and Schramm \cite{BS}. For a vertex $p\in V(G_k)$ we consider the
rooted $K$-colored Schreier graph $(G_k^{c_k},p)$. The measure $\mu_k$
is defined as
$$\mu_k=\frac{1}{|V(G_k)|} \sum_{p\in V(G_k)} \delta(G_k^{c_k},p)\,,$$
where $\delta(G_k^{c_k},p)$ is the Dirac-measure on $\Sch^{\overline{Q},K}_{\F_n}$
concentrated on the rooted $K$-colored Schreier graph $(G_k^{c_k},p)$.
Clearly, $\mu_k$ is invariant under the action of $\F_n$.
Since the space of $\F_n$-invariant probability measures on the
compact space $\Sch^{\overline{Q},K}_{\F_n}$ is compact with respect to
the weak-topology, we have a convergent subsequence $\{\mu_{n_k}\}^\infty_{k=1}$
converging weakly to some probability measure $\mu$.
We consider the $K$-Bernoulli shift space $B^K(Z)$ as an $\F_n$-space, where for $h\in\F_n$ and  $f\in B^K(Z)$
$$\overline{\beta}(h)(f)=\beta(\kappa(h))(f)\,,$$
\noindent
where $\beta$ is the left $\G$-action on $B^K(Z)$.
Hence, we have an injective $\G$-equivariant map $\Phi_\kappa:B^K(Z)\to \Sch^{\overline{Q},K}_{\F_n}$.
\begin{lemma}\label{conc1}
The probability measure $\mu$ is concentrated on the $\F_n$-invariant
closed set $\Omega$ of nonrepetitive $K$-colorings in $\Phi_\kappa(B^K(Z))$.
\end{lemma}
\proof
Let $U_r\subset \Sch^{\overline{Q},K}_{\F_n}$ be the clopen set of
$K$-colored Schreier graphs $G$ such that
the ball $B_r(G,\roo(G))$ is not rooted-labeled isomorphic to 
$B_r(S^Q_\G(H),H)$ for some $H\in Z$.
By our assumptions on the sofic approximations, $\lim_{k\to\infty}\mu_k(U_r)=0\,,$
hence $\mu(U_r)=0\,.$
Now let $V_r\subset \Sch^{\overline{Q},K}_{\F_n}$ be the clopen set of
$K$-colored Schreier graphs $G$ such that the ball
$B_r(G,\roo(G))$ contains a repetitive path. By our assumptions on the
colorings $c_k$, $\mu_k(V_r)=0$ for any $k\geq 1$. Hence $\mu(V_r)=0$.
Therefore $\mu$ is concentrated on $\Omega$. \qed

\vi
Now we can finish the proof of Theorem \ref{tetel2}.
We can identify $\Omega$ with a $\G$-invariant closed subset $\overline{\Omega}$ of $B^K(Z)$ on
which the $\G$-action is $Z$-proper by Proposition \ref{propo2}. That is, our construction gave rise
to a $Z$-proper Bernoulli subshift with an invariant measure. \qed

\noindent
Note that we have a $\G$-equivariant continuous map from the $Z$-proper space above to $Z$ itself
mapping $x$ into $\Stab(x)$. Recall that a $\G$-invariant
measure on $\Sub(\Gamma)$ is called an invariant random subgroup.
\begin{propo}
Any sofic URS admits an invariant measure.
\end{propo}

\noindent
Recall that a $\G$-invariant
measure on $\Sub(\Gamma)$ is called an invariant random subgroup \cite{AGV}. Example 3.3 in \cite{GW} shows 
that there exists a uniformly recurrent subgroups $Z\subset \Sub(\F_2)$ that does not admit invariant random subgroups, hence $Z$ is not sofic. In Section \ref{example}, we provide further examples of uniformly recurrent
subgroups that does not carry invariant measures. 
\section{Coamenable uniformly recurrent subgroups} \label{example}
\subsection{Colored graphs} \label{colored}
Let $\G_k$ be the $k$-fold free product of cyclic groups
of rank $2$, with free generators $A=\{a_i\}^k_{i=1}$. 
Let $G$ be an arbitrary infinite, simple connected graph of bounded vertex degrees and a proper edge-coloring by
$k$-colors $\{c_1,c_2,\dots,c_k\}$. Observe that the edge-coloring
of $G$ (and picking an arbitrary root) gives rise to a $\G_k$-Schreier graph
$(S,x)$. 
 The action of $\G_k$ on $V(G)$ is defined the following way. If $x\in V(G)$ and $1\leq i\leq k$, then
\begin{itemize}
\item If there is no $c_i$-colored edge adjacent to $x$, then $a_i(x)=x$.
\item If there there exists an edge $(x,y)$ colored by $c_i$, then $a_i(x)=y$.
\end{itemize}

\noindent
Let $X$ be the orbit closure of the rooted Schreier graph $(S,x)$ above. Then it contains a minimal system $(M,\G_k,\beta)$.
Then $(S^Q_{\G_k})^{-1}(M)$ is a uniformly recurrent subgroup, where $S^Q_{\G_k}:\Sub(\G_k)\to
\Sch^Q_{\G_k}$ is the map defined in Subsection \ref{rooted}.
\begin{propo}\label{constru}
For any infinite simple, connected graph $G$ of bounded vertex degree, there exists $k>0$ and a edge-coloring of $G$
with $k$ colors such that all the uniformly recurrent subgroups that can be obtained as above are necessarily generic.
\end{propo}
\proof
First, consider an arbitrary proper edge-coloring $c:E(G)\to L$ and a proper nonrepetitive vertex-coloring $\rho:V(G)\to D$ by some finite
sets $L$ and $D$ (the product of a nonrepetitive and a proper vertex-coloring
is always a proper nonrepetitive vertex-coloring). Now we construct a new proper edge-coloring $\zeta$ of $G$ by the set $D_2\times L$, where
$D_2$ is the set of $2$-elements subset of $D$.
Let $\zeta(e)=\{\rho(x),\rho(y)\}\times c(e)$,
where $x,y$ are the endpoints of $e$. Since $\rho$ is proper, $\rho(x)\neq \rho(y)$.
Hence we obtain a Schreier graph $T\in \Sch^A_{\G_k}$, where $k=|D_2\times L|$ and
$A=\{a_1,a_2,\dots,a_k\}.$
Let $(M,\G_k,\beta)$ be a minimal subsystem in the orbit closure of $T$. By Proposition \ref{assy}, it
is enough to show that if $T'\in M$ and $x\neq y\in V(T')$, then $(T',x)$ and $(T',y)$ are not
rooted-labeled isomorphic.
We construct a nonrepetitive vertex-coloring $\rho':V(T')\to D$ in the following way.
If $\deg(z)>1$, let $\rho'(z)=d$, where $d$ is the unique element in the intersection of the $D_2$-components
of the edges adjacent to $z$. If $\deg(z)=1$, then let $\rho'(z)=d$, where the $D_2$-component
of $z$ is $\{d,d'\}$ and $\rho'(z')=d'$, for the only neighbour of $z$. Since $T'$ is in the orbit closure
of $T$, the coloring $\rho'$ is nonrepetitive, hence by Proposition \ref{propo2}, $T'(x)$ and $T'(y)$ are not rooted-labeled isomorphic. \qed.

\subsection{Coamenability} \label{coame1}
Let $\G$ be a finitely generated group and $H\in\Sub(\G)$. Recall that $H$ is {\bf coamenable} if the action of $\G$
on $\G/H$ is amenable. That is, there exists a sequence of finite subsets $\{F_k\}^\infty_{k=1}\subset \G/H$
such that for any $g\in\G$,
$$\lim_{n\to \infty} \frac{|gF_k\cup F_k|}{|F_k|}=1\,.$$
\noindent
We call a URS $Z\subset\Sub(\G)$ coamenable if for all $H\in Z$, $H$ is coamenable.
\begin{propo}
Let $Z\subset\Sub(\G)$ be a URS such that there exists $H\in Z$ so that $H$ is coamenable.
Then $Z$ is coamenable.
\end{propo}
\proof Fix a generating system $Q=\{\gamma_i\}^n_{i=1}$ for $\G$.
Let $\{F_k\}^\infty_{k=1}$ be finite subsets in $\G/H$ such that for any $g\in\G$, 
$\lim_{n\to \infty} \frac{|gF_k\cup F_k|}{|F_k|}=1\,.$
Let $x_k\in F_k$ and $l(k)>0$ such that if $y\in F_k$ then $B_{l(k)}(\G/H,x_k)$ contains
$B_k(\G/H,y)\,.$ Let $K\in Z$. Since $Z$ is a URS, for any $k\geq 1$, there
exists $x'_k\in \G/K$ such that
$B_{l(k)}(\G/H,x_k)$ is rooted-labeled isomorphic to $B_{l(k)}(\G/K,x'_k)$. For $k\geq 1$, let $F_k'\subset \G/K$ be the image
of $F_k$ by the isomorphism above.
If $g\in\G$ then $gF_k'\subset B_{l(k)}(\G/K,x_k')$ provided that $k$ is large enough (depending on $g$). Hence,
$\lim_{k\to \infty} \frac{|gF'_k\cup F'_k|}{|F'_k|}=1\,.$  \qed
\vskip 0.2in
\noindent
Let $G$ be an arbitrary graph that is amenable in the sense that there exists a sequence of subsets
$\{F_k\}^\infty_{k=1}$ such that $\lim_{k\to\infty} \frac{|B_1(F_k)|}{|F_k|}=1\,,$ where $x\in B_1(F_k)$ if either $x\in F_k$
or there exists $y\in F_k$ adjacent to $x$. Repeating the proof of the previous proposition one can immeadiately
see that all the $\G_k$-URS's constructed from $G$ as in Subsection \ref{colored} are coamenable.
Barlow [Proposition 4.\cite{Barl}] has shown that for any $\alpha\geq 1$ there exists a bounded degree infinite graph $G_\alpha$ and
 positive constants
$C^1_\alpha$ and $C^2_\alpha$ such that 
\begin{equation}
C^1_\alpha r^\alpha \leq B_r(G_\alpha,x) \leq C^2_\alpha r^\alpha
\end{equation}
\noindent
holds for all $x\in V(G_\alpha)$. Such graphs are clearly amenable. Hence we can see that as opposed to 
finitely generated group case, for any $\alpha\geq 1$ there exist generic coamenable uniformly recurrent subgroups $Z\subset\Sub(\G_k)$ so that the volume growth rate of the individual Schreier graphs $S^A_{\G_k}(H), H\in Z$ are always $\alpha$.
\subsection{Coamenable uniformly recurrent subgroups are sofic}\label{sub43}
The following theorem (or rather the construction in the proof) will be crucial in Section \ref{amena}.
\begin{theorem}
Let $\G$ be a finitely generated group and $Z\subset\Sub(\G)$ be a coamenable URS. Then $Z$ is sofic.
\end{theorem}
\proof Fix a generating system $Q=\{\gamma_i\}^n_{i=1}$. Again, let $\F_n$ be the free group with
free generating system $\overline{Q}=\{r_i\}^n_{i=1}$ and $\kappa:\F_n\to\Gamma$ be the corresponding
quotient map. Every continuous action $\alpha$ of $\Gamma$ can be regarded as a $\F_n$-action $\alpha\circ\kappa$.
In particular, we have a $\F_n$-invariant embedding $\lambda:S^Q_\G(Z)\to \Sch^{\overline{Q}}_{\F_n}\,.$ Now let $H\in Z$
and consider the Schreier graph $S^Q_\G(H)$.
Since $Z$ is coamenable, the isoperimetric constant of $S^Q_\G(H)$ is zero, that is, we have a sequence
of finite induced subgraphs $\{H_k\}^\infty_{k=1}\subset S^Q_\G(H)$ so that
$$\lim_{k\to\infty} \frac{|\partial H_k|}{|V(H_k)|}=0\,,$$
\noindent
where $\partial H_k$ is the set of vertices in $H_k$ for which there exists \\$y\in V(S^Q_\G(H))\backslash V(H_k)$
with $\gamma_ix=y$ or $\gamma_i y=x$ for some $1\leq i\leq n$.
Now we construct a sequence of $\F_n$-Schreier graphs $\{G_k\}^\infty_{k=1}$ that form a sofic approximation
of $Z$:
\begin{itemize}
\item $V(G_k)=V(H_k)$.
\item If $\gamma_i x=y$ for some $x,y\in V(H_k)$, then $r_ix=y$.
\item Then the action of $r_i$ is extended to the set $V(G_k)$ arbitrarily.
\end{itemize}
\noindent
Let $W_r(G_k)\subset V(G_k)=V(H_k)$ be the set of vertices $p$ for which $d_{G_k}(p,\partial H_k)>r$. Clearly, if $p\in W_r(G_k)$
then $B_r(W(G_k),p)$ is rooted-labeled isomorphic to  $B_r(S^Q_\Gamma(H),p)$. That is, all the vertices of $W_r(G_K)$ are $(Z,r)$-vertices.
Since, $\frac{|\partial H_k|}{|V(G_k)|}\to 0$, we have that $\frac{|W_r(G_k)|}{|V(G_k)|}=1$. Hence the Schreier graphs $\{G_k\}^\infty_{k=1}$ form a sofic approximation
of the URS Z. \qed
\vskip 0.1in
\noindent
As in Subsection \ref{sub42}, for each $k\geq 1$ we have an $\F_n$-invariant probability measure on
$\Sch^{\overline{Q}}_{\F_n}$
$$\mu_k=\frac{1}{|V(G_k)|}\sum_{p\in V(G_k)} \delta(G_k,p)\,.$$
\noindent
Let $\{\mu_{n_k}\}^\infty_{l=1}$ be a weakly convergent sequence converging to an $\F_n$-invariant measure $\mu$
on $\Sch^{\overline{Q}}_{\F_n}$. By our previous lemma, the measure $\mu$ is concentrated on $\lambda(S^Q_\G(Z))$.
Note that the probability measure $\mu$ depends only on the sequence of subgraphs $\{H_{k_l}\}^\infty_{l=1}$.
We say that the sequence of subgraphs $\{H_l\}^\infty_{l=1}$ is {\bf
  convergent in the sense of Benjamini and Schramm} if the associated
probability measures $\{\mu_l\}^\infty_{l=1}$ converge to some invariant measure $\mu$ on $\lambda(S^Q_\G(Z))$.
In this case the measure preserving action $(Z,\G,\lambda^{-1}(\mu))$ is called the {\bf limit} of the sequence  
$\{H_l\}^\infty_{l=1}$. Also note, that if $Z$ is a generic URS, then  $(Z,\G,\lambda^{-1}(\mu))$ is a totally nonfree action
in the sense of Vershik \cite{Vers}.
\subsection{A characterization of coamenability}
As in the previous sections let $\G$ be a finitely generated
group with generating system $Q=\{\gamma_i\}^r_{i=1}$
and $Z\subset\Sub(\Gamma)$ be a URS.
The goal of this subsection is to prove the following characterization of
coamenability.
\begin{theorem}\label{karak}
The URS $Z$ is coamenable if and only if
every $Z$-proper continuous action of $\G$ admits
an invariant measure.
\end{theorem}
\proof
First, let $Z$ be coamenable and $\alpha:\G\curvearrowright X$ be a continuous
$Z$-proper
action. Let $x \in X$ and $Stab_\alpha(x)=H\in Z\,.$
Then the orbit graph of $x$ is isomorphic to $S^Q_\G(H)$.
Let $\{F_k\}^\infty_{k=1}\subset \G/H$ be a sequence of 
finite sets such that for any $g\in\G$, 
\begin{equation}\label{folni}
\lim_{k\to\infty}
\frac{|gF_k\cup F_k|}{|F_k|}=1\,.
\end{equation}
Now we proceed in exactly the same way as in the proof of the classical
Krylov-Bogoliubov Theorem.  Fix an ultrafilter $\omega$ on the natural numbers
and let $\lim_\omega$ be the corresponding ultralimit. We define a bounded
linear
functional $T:\C[X]\to \C$ in the following way. For a continuous function
$f:X\to \C$ set 
$$T(f)=\lim_\omega \frac{\sum_{\overline{\gamma}\in F_k}
  f(\overline{\gamma}x)}{|F_k|}\,.$$
\noindent
Then, by (\ref{folni}), $T(\gamma(f))=T(f)$ for all $\gamma\in\G$, $T(1)=1$,
therefore $T:\C[X]\to \C$ is a $\G$-invariant bounded functional associated to
a $\G$-invariant Borel probability measure $\mu$.
\noindent
Now, let $Z$ be a URS that is not coamenable and let $H\in Z$.
Then the graph $S^Q_\G(H)$ has positive isoperimetric constant so by Theorem 3.1
\cite{Sos},
we have maps $\phi_1:\Gamma/H\to \Gamma/H$, $\phi_2:\Gamma/H\to \Gamma/H$
so that $\phi_1(\Gamma/H)\cap \phi_2(\Gamma/H)=\emptyset$ and there exists a
positive constant $C>0$ so that for all $p\in\G/H$
$$d_{S^Q_\G(H)} (\phi_1(p),p))< C\quad \mbox{and} \quad d_{S^Q_\G(H)}
(\phi_2(p),p))<C\,.$$
\noindent
Now we build a vertex-coloring for the graph $S^Q_\G(H)$ to encode $\phi_1$ and
$\phi_2$.
First, we pick a nonrepetitive coloring $c_1:\G/H\to D$, where $D$ is some
finite set. Then we choose a coloring $c_2:\G/H\to E$ for some finite set $E$
so that $c_2(p)\neq c_2(q)$, whenever $$0< d_{S^Q_\G(H)}(p,q)\leq 3C\,.$$
\noindent
We need two more colorings of the vertices of $S^Q_\G(H)$:
$$c_3:\G/H\to \{E\times \{1\}\}\cup \{*\}$$
and
$$c_4:\G/H\to \{E\times \{2\}\}\cup \{*\}$$
\noindent
satisfying the following properties.
If there exists $p\in \G/H$ so that $\phi_1(p)=q$ and $c_2(p)=e$, then
$c_3(q)=e\times\{1\}$. If such $p$ does not exist, set 
$c_3(q)=*$.
If there exists $p\in \G/H$ so that $\phi_2(p)=q$ and $c_2(p)=e$, then
$c_4(q)=e\times\{2\}$. If such $p$ does not exist, set 
$c_4(q)=*$.
\noindent
Let $M=D\times E\times  \{\{E\times \{1\}\}\cup \{*\} \}\times \{ \{E\times
\{2\}\}\cup \{*\}\}$.
Our final coloring $c:\G/H\to M$ is defined by 
$$c(p)=c_1(p)\times c_2(p)\times c_3(p)\times c_4(p)\,.$$
\noindent
Let $X$ be the orbit closure of the $M$-colored graph $S_\G^{Q,c}(H)$ in
the space $\Sch^{M,Q}_\G\,.$
Observe that $c$ is nonrepetitive since even its first component is
nonrepetitive, that is the action $\beta$ on $X$ is $Z$-proper.
Now we need to show that $\beta$ admits no $\G$-invariant measure.
We define continuous injective maps $\Phi_1:X\to X$ and $\Phi_2:X\to X$ such 
that
\begin{itemize}
\item $\Phi_1(X)\cap \Phi_2(X)=\emptyset$
\item For each $x\in X$, $\Phi_1(x),\Phi_2(x)\in \Orb(x)\,.$
\end{itemize}
\noindent
Thus the equivalence relation defined by the action is compressible, so it
cannot admit an invariant measure \cite{DJK}. 
The construction of $\Phi_1$ and $\Phi_2$ goes as follows.
If $x\in X$, then $c(x)=(c_1(x),c_2(x),c_3(x),c_4(x))$ is well-defined and
there exist a unique $y\in X$ and a unique $z \in X$ such that
\begin{itemize}
\item $c_2(x)\times {1}=c_3(y)$, $c_2(x)\times {2}=c_4(z)$
\item $d_{\Orb(x)}(x,y)\leq C$, $d_{\Orb(x)}(x,z)\leq C$.
\end{itemize}
\noindent
We set $\Phi_1(x)=y$, $\Phi_2(x)=z$, Clearly, $\Phi_1$ and $\Phi_2$ are
continuous and $\Phi_1(X)\cap \Phi_2(X)=\emptyset\,.$ \qed
\vskip 0.2in
\noindent
Let $\G$ be as above and $Z\subset \Sub(\G)$ be a URS that is not coamenable
and $X$ is a $Z$-proper action without invariant measure as above.
Let $Y\subset X$ be a minimal $Z$-proper $M$-Bernoulli subshift.
Let $\G_{M}$ be the $M$-fold free product of the finite group of two elements
with free generator system $\{a_m\}_{m\in M}$.
Then we can associate to $Y$ a nonsofic generic URS $Z\subset \Sub(\G * \G_M)$  
in the following way.
Let $S=S_\G^{M,Q}(H)$ be an element of $Y$. Let $V=V(S)\times \{0,1\}$.
We define an action of the group $\G * \G_M$ as follows.
The group $\G$ acts on $V(S)\star \{0\}$ as $\G$ acts on $V(S)$. Also, the
group $\G$ acts on $V(S)\star \{1\}$ trivially.  If
$x\in V(S)$, $c(x)=m$, then $a_m(x\times\{0\})=x\times\{1\}$ and
$a_m(x\times\{1\})=x\times\{0\}$.
Otherwise, let $a_m(x\times \{0\})=x\times\{0\}$ and $a_m(x\times
\{1\})=x\times\{1\}$.
It is not hard to see that the resulting $(\G * \G_M)$-Schreier graph satisfies the 
conditions of Proposition \ref{repet} and \ref{assy}, hence the associated URS
is generic and does not admit invariant measures.

\section{The $C^*$-algebras of uniformly recurrent subgroups}
\subsection{The algebra of local kernels}
Let $\G$ be a finitely generated group with generating system $Q=\{\gamma_i\}^n_{i=1}$
and $Z\subset\Sub(\Gamma)$ be a URS of $\G$. Let $H\in Z$ and $S=S^Q_\G(H)$ be the Schreier graph
of $H$. A {\bf local kernel} is a function $K:\G/H \times \G/H\to \C$
satisfying the following properties.
\begin{itemize}
\item There exists an integer $R>0$ (depending on $K$) such that
$K(x,y)=0$ if $d_S(x,y)>R$.
\item If $B_R(S,y)$ is rooted-labeled isomorphic to $B_R(S,z)$ then
$K(y,\gamma y)=K(z,\gamma z)$ provided that $d_S(y,\gamma y)=d_S( z, \gamma z)\leq R$.
\end{itemize}
\noindent We will call the smallest $R$ satisfying the two conditions above the {\bf width} of $K$.
\noindent
It is easy to see that the local kernels form a unital $*$-algebra $\C Z$
with respect to the following operations:
\begin{itemize}
\item $(K+L)(x,y)=K(x,y)+L(x,y)$
\item $KL(x,y)=\sum_{z\in \G/H} K(x,z)L(z,y)$ 
\item $K^*(x,y)=\overline{K(y,x)}$.
\end{itemize}
\noindent
By minimality, the algebra $\C Z$ does not depend on the choice of
$H$ or the generating system $Q$ only on the URS $Z$ itself. We will call the concrete realization
of the algebra of local kernels az above the representation of $\C Z$ on $\C^{\G/H}$.
One can observe that if $Z$ consists only of the unit element, then $\C Z$ is the complex group algebra
of $\Gamma$.
\subsection{The construction of $C^*_r(Z)$}
Let $\G$ be a finitely generated group (with a fixed generating system $Q=\{\gamma_i\}^n_{i=1}$) and $Z\subset\Sub(\G)$ be
a URS. Let $H\in Z$ and consider the algebra $\C Z$ as above represented on the vector space $\C^{\G/H}$. The we have
a bounder linear representation of $\C Z$ on $l^2(\G/H)$ by
$$K(f)(x)=\sum K(x,y)f(y)\,,$$
\noindent
where $f\in l^2(\G/H)$. 
\begin{defin}
The $C^*$-algebra of $Z$, $C^*_r(Z)$  is defined as the norm closure of $\C Z$ in $B(l^2(\G/H))$.
\end{defin}
\noindent
Note that we used a specific subgroup $H$ in order to equip the algebra $\C Z$ with a norm. However, we have the following
proposition.
\begin{propo}
The norm on $\C Z$ and hence the definition of  $C^*_r(Z)$ does not depend on the choice of the subgroup $H$.
\end{propo}
\proof
Let $K\in \C Z$ be a local kernel of width $R$ and let $H,L\in Z$. Let $K_H$ respectively $K_L$ be the
representation of $K$ on $l^2(\G/H)$ respectively on $l^2(\G/L)$. We need to show that
$\|K_H\|=\|K_L\|$\,. Let $\e>0$ and $f\in l^2(\G/H), \|f\|=1$ such that
$f$ is supported on a ball $B_T(S^Q_\G(H),x)$ and
$\|K_H(f)\|\geq \|K_H\|-\e\,.$ Observe that
$K_H(f)$ is supported on the ball $B_{T+R}(S^Q_\G(H),x)$ and
$\|K_H(f)\|\geq \|K_H\|-\e\,.$ By Proposition \ref{repet}, there exists $y\in\G/L$ such that
the balls $B_{T+R}(S^Q_\G(H),x)$ and
$B_{T+R}(S^Q_\G(L),y)$ are rooted-labeled isomorphic. Hence, there exists $f'\in l^2(\G/L)$ supported
on $B_T(S^Q_\G(L),y)$, $\|f'\|=1$ such that $\|K_H(f)\|=\|K_L(f')\|\,.$ Therefore, $\|K_H\|\leq \|K_L\|$.
Similarly, $\|K_L\|\leq \|K_H\|$, that is, $\|K_H\|= \|K_L\|$.\qed
\subsection{The $C^*$-algebras of generic URS's are simple}\label{s63}
The goal of this section is to prove the following theorem.
\begin{theorem}\label{ursimple}
Let $\G$ be as above and $Z\subset\Sub(\G)$ be a generic URS. Then the $C^*$-algebra $C^*_r(Z)$ is simple. 
\end{theorem}
\proof
Let $H\in Z$. For each $r\geq 1$ we define an equivalence relation on $\G/H$ in the
following way. If $p,q\in\G/H$, then $p\equiv_r q$ if
the balls $B_r(S^Q_\G(H),p)$ and $B_r(S^Q_\G(H),q)$ are rooted-labeled isomorphic.
The following lemma is a straightforward consequence of Proposition \ref{repet} and Proposition \ref{assy}.
\begin{lemma} \label{alapok}
Let $\equiv_r$ be the equivalence relation as above. Then:
\begin{enumerate}
\item For any $n\geq 1$ there exists $r_n$ such that if $p\neq q$ and $p\equiv_{r_n} q$, then
$d_{S^Q_{\G}(H)} (p,q )\geq n\,.$
\item For every $r\geq 1$ there exists $t_r$ such that for any $p\in \G/H$ the ball
$B_{t_r}(S^Q_\G(H),p)$ intersects all the equivalence classes of $E_r$ (in particular, the number of equivalence classes
is finite).
\item If $r\leq s$, then $p\equiv_s q$ implies $p\equiv_r q$.
\item Let $E_r$ denote the classes of $\equiv_r$. Then we have an inverse system of surjective maps
$$E_1\leftarrow E_2\leftarrow \dots$$
and a natural homeomorphism $\iota_H:\lim_{\leftarrow} E_r\to Z$, between the
compact space $\lim_{\leftarrow} E_r$ and the uniformly recurrent subgroup $Z$.
\end{enumerate}
\end{lemma}
\noindent
Note that if $\alpha\in E_r$, then $\iota_H(\alpha)$ is the clopen set of
Schreier graphs $S^Q_\Gamma(L)$, $L\in Z$, such that the ball $B_r(S^Q_\G(L),L)$ is
rooted-labeled isomorphic to the ball $B_r(S^Q_\G(H),x)$, where
$x\in\alpha\,.$

\noindent
Now let us consider the commutative $C^*$-algebra $l^\infty(\G/H)$. For any $r\geq 1$ and $\alpha\in E_r$ we
have a projection $e_\alpha\in l^\infty(\G/H)$, where $e_\alpha(x)=1$ if $x\in\alpha$ and zero otherwise. The 
projections $\{e_\alpha\}_{r\geq 1,\alpha\in E_r}$ generates a *-subalgebra $\cA$ in $l^\infty(\G/H)$ and by the
previous lemma the closure of $\cA$ in $l^\infty(\G/H)$ is isomorphic to $\C[Z]$ (the $C^*$-algebra
of continuous complex-valued functions on the compact metrizable space
$Z$). Indeed, under this isomorphism $\lambda_H:\overline{\cA}\to\C[Z]$,
$\lambda_H(e_\alpha)$ is the
characteristic function of the clopen set $\iota_H(\alpha)$.
It is easy to see that the isomorphism $\lambda_H:\overline{\cA}\to \C[Z]$ commutes with the respective
$\G$-actions.
Now let us consider the representation of $C^*_r(Z)$ on $l^2(\G/H)$. For $K\in C^*_r(Z)$ let
$K(x,y)=\langle K(\delta_y),\delta_x\rangle$, be the kernel of $K$.
We have a bounded linear map $Q_r: C^*_r(Z)\to C^*_r(Z)$ given by
$$Q_r (K)=\sum_{\alpha\in E_r} e_\alpha K e_\alpha\,.$$
\begin{lemma}
For any $r\geq 1$, $\| Q_r\|\leq 1\,.$
\end{lemma}
\proof
Let $h\in l^2(\Gamma/H)$, $\|h\|=1$. For any $K\in C^*_r(Z)$ we have
that
$$\|(Q_r(K))(h)\|^2=\| \sum_{\alpha\in E_r} e_\alpha K e_\alpha (h)\|^2=
\sum_{\alpha\in E_r} \|e_\alpha K e_\alpha (h)\|^2\leq $$
$$\leq \|K\|^2 \sum_{\alpha\in E_r} \|e_\alpha(h)\|^2=\|K\|^2\,.$$
\noindent
Therefore $\|Q_rK\|\leq \|K\|\,.$ \qed
\vskip 0.1in
\noindent
Observe that we have a natural injective homomorphism
$\rho:\cA\to\C Z$ defined in the following way.
\begin{itemize}
\item $\rho(a)(x,x)=a(x)\,.$
\item $\rho(a) (x,y)=0$ if $x\neq y$.
\end{itemize}
\noindent
Clearly, $\rho$ is preserving the norm, so we can extend it to a unital embedding
$\overline{\rho}:\cA\to C^*_r(Z)\,.$
Also, we have a map $\kappa:\G\to C^*_r(Z)$ such that
$\kappa(g)(x,y)=1$, whenever $g^{-1}x=y$ and $\kappa(g)(x,y)=0$ otherwise.
\begin{lemma}\label{kappa}
For any $g,h\in\G$, $\kappa(g)\kappa(h)=\kappa(gh)\,.$
\end{lemma}
First, we have that
$$\kappa(g)\kappa(h)(x,y)=\sum_{z\in\G/H} \kappa(g)(x,z)\kappa(h)(z,y)\,.$$
\noindent
Hence, $\kappa(g)\kappa(h)(x,y)=1$ if $y=h^{-1}g^{-1}x$ and $\kappa(g)\kappa(h)(x,y)=0$ otherwise.
Therefore, $\kappa(g)\kappa(h)=\kappa(gh)\,.$ \qed
\begin{lemma} \label{equiv}
For any $g\in\G$ and $a\in\overline{\cA}$
$$\rho(g(a))=\kappa(g)\rho(a)\kappa(g^{-1})\,.$$
\end{lemma}
\proof
On one hand,
$\rho(g(a))(x,y)=a(g^{-1}(x))$ if $x\neq y$, \\ otherwise $\rho(g(a))(x,y)=0\,.$
On the other hand,
$$\kappa(g)\rho(a) \kappa^{-1}(g)(x,x)=\sum_{y\in\Gamma/H} \kappa(g)(x,y)\rho(a)(y,y)
\kappa^{-1}(g)(y,x)=a(g^{-1}(x))\,.$$
Also, $\kappa(g)\rho(a) \kappa^{-1}(g)(x,y)=0$ if $x\neq y$. \qed
\vskip 0.2in
\noindent
Let us consider the linear operator $D:C^*_r(Z)\to C^*_r(Z)$ such that
for $x\in\G/H$ $D(K)(x,x)=K(x,x)$, $D(K)(x,y)=0$ if $x\neq y$. The operator $D$ is bounded with
norm $1$
since 
$$\|D(K)\|=\sup_{x\in\G/H} |K(x,x)|=\sup_{x\in\G/H}|\langle K(\delta_x),\delta_x\rangle|\,.$$
\begin{lemma}\label{limd}
Let $K\in \C Z$. Then $Q_r(K)=D(K)$ provided that $r$ is large enough.
\end{lemma}
\proof
Let $s>0$ be the width of $K$ and let $r>0$ be so large that if $p\equiv_r q$ and $p\neq q$, then
$d_{S^Q_\Gamma(H)}(p,q)>s$. Then, if $\alpha\in E_r$ we have that $(e_\alpha K e_\alpha)(x,y)=0$ if $x\neq y$ or $x\notin 
\alpha$, otherwise $(e_\alpha K e_\alpha)(x,x)=K(x,x)$. Therefore, $Q_r(K)=D(K).$ \qed
\begin{lemma}
Let $K\in C^*_r(Z)$. Then $\lim_{r\to\infty} Q_r(K)=D(K)\,.$
\end{lemma}
\proof
Let $K_n\to K$ such that $K_n\in\C Z$. Then, by the previous lemma
$\|Q_r(K)-D(K_n)\|\leq \|K-K_n\|$, provided that $r$ is large enough.
Since $D(K_n)\to D(K)$, we have that $\lim_{r\to\infty} Q_r(K)= D(K)\,.$ \qed
\begin{lemma}\label{utolso}
Let $I \lhd C^*_r(K)$ be a closed ideal. Suppose that $I\cap D(C^*_r(Z))\neq \{0\}$. Then
$I=C^*_r(K)$.
\end{lemma}
\proof
Recall that $D(C^*_r(Z))=\rho(\overline{\cA})$, so by Lemma \ref{equiv} we have
a nonzero, $\G$-invariant closed ideal in $\cA\cong \C[Z]$.
However, any $\G$-invariant closed ideal in $\cA\cong \C[Z]$ is in the form of $I(Y)$,
where $Y$ is a $\G$-invariant closed set in $Z$ and $I(Y)$ is the set of continuous functions
vanishing at $Y$. By minimality, $Y$ must be empty, hence $I$ contains the unit, that is, 
$I=C^*_r(K)$. \qed
\vskip 0.2in
\noindent
Now, we finish the proof of our theorem.
Let $I$ be a closed ideal of $C^*_r(Z)$ and $0\neq K\in I$. Then $K^*K\in I$ and
$D(K^*K)\neq 0$. Since $D(K^*K)=\lim_{r\to \infty} Q_r(K^*K)$ and $Q_r(K^*K)\in I$ for any $r\geq 1$,
we have that $1\in I$. \qed

\vskip 0.2in
\noindent
{\bf Remark} Let $Z\subset\Sub(\Gamma)$ be a not necessarily generic URS, where
$\G$ is a finitely generated group as above. Let
$Y\subset S^{K,Q}_\G(Z)$ be a minimal $Z$-proper Bernoulli subshift. Then the
local kernels on $Y$ can be defined using the rooted-labeled-colored neighborhoods
and the resulting $C^*$-algebra is always simple. 

\section{Exactness and nuclearity}
\subsection{Property A vs. Local Property A}
First let us recall the notion of {\bf Property $A$} from \cite{Roe}.
Let $G$ be an infinite graph of bounded vertex degrees. We say the $G$ has
Property A if there exists a sequence of maps $\{\vs^n:V(G)\to l^2(V(G)\}^\infty_{n=1}$
such that
\begin{itemize}
\item Each $\vs^n_x$ has length $1$.
\item If $d_G(x,y)\leq n$, then $\|\vs^n_x-\vs^n_y\|\leq \frac{1}{n}\,.$
\item For any $n\geq 1$ we have $R_n>0$ such that the vector $\vs^n_x$ is
  supported in the ball $B_{R_n}(G,x)$.
\end{itemize}
\vskip 0.2in
We also need the notion of the {\bf uniform Roe algebra} of the graph $G$.
First, we consider the $*$-algebra of {\bf bounded kernels} $K:V(G)\times
V(G)
\to\C$, that is
\begin{itemize}
\item there exists some positive integer $R$ depending on $K$ such that $K(x,y)=0$ if
  $d_G(x,y)>R$,
\item there exists some positive integer $M$ depending on $K$ such that \\
  $|K(x,y)|<M$.
\end{itemize}
\noindent
The uniform Roe algebra $C^*_u(G)$ is the norm closure of the bounded kernels in
$B(l^2(V(G)))$. Observe that if $Z\subset\Sub(\Gamma)$ is a unformly recurrent
subgroup and $H\in Z$, $S=S^Q_\G(H)$, then $C^*_r(Z)\subset C^*_u(S)$.
According to Proposition 11.41 \cite{Roe}, if $G$ has Property $A$ then the
algebra $C^*_u(G)$ is nuclear.
All $C^*$-subalgebras of a nuclear $C^*$-algebra are exact, hence we have the
following proposition.
\begin{propo}\label{maexact}
Let $Z\subset\Sub(\Gamma)$ and $H\in Z$ as above, such that $S^Q_\G$ has Property
$A$.
Then $C^*_r(Z)$ is exact.
\end{propo}
\vskip 0.2in
\noindent
{\bf Example:} Let $G$ be the underlying graph of the Cayley graph of an exact group (say, a
hyperbolic group or an amenable group) and let $S$ be a colored graph
associated to a generic URS $Z$ as in Proposition \ref{constru}. Then by the
previous proposition, $C_r(Z)$ is a simple exact $C^*$-algebra.
\vskip 0.2in
\noindent
Now we introduce the notion of Schreier graphs with {\bf Local Property $A$}. 
\begin{defin}
Let $S=S^Q_\G(H)$ be a Schreier graph. We say that $S$ has Local Property
$A$, if the sequence $\{\vs^n\}^\infty_{n=1}$ can be chosen locally, that is
for any $n\geq 1$, there exists $S_n> R_n$ so that for $x,y\in V(G)$ the balls
$B_{S_n}(G,x)$ and $B_{S_n}(G,y)$ are rooted-labeled isomorphic under the map
$\theta:B_{S_n}(G,x)\to B_{S_n}(G,y)$, then $\vs^n_y=\theta(\vs^n_x)$.
\end{defin}
\noindent
The
main result of this section is the following theorem.
\begin{theorem} \label{manuclear}
Let $\G$ be a finitely generated group, $Z\subset\Sub(\G)$ a uniformly recurrent subgroup and $H\in Z$ so that
$S^Q_{\G}(H)$ has local Property A. Then $C^*_r(Z)$ is nuclear.
\end{theorem}
\proof
We closely follow the proof of Proposition 11.41 \cite{Roe}. The nuclearity
of the uniform Roe algebra for a graph $S$ having Property $A$ has been proved
the following way (we will denote by $X$ the vertex set of $S$).
First, a sequence of unital completely positive maps $\Phi_n:C^*_u(S)\to
l^\infty(X)\otimes M_{N_n}(\C)$ were constructed, where $M_{N_n}(\C)$ is the algebra of
$N_n\times N_n$-matrices.
Then, a sequence of unital completely positive maps $\Psi_n:l^\infty(X)\otimes
M_{N_n}\to C^*_u(S)$ were given in such a way that
$\{\Psi_n\circ\Phi_n\}^\infty_{n=1}$ tends to the identity in the point-norm
topology. Hence, the nuclearity of the uniform Roe algebra $C^*_u(S)$ follows.
It is enough to see that $\Phi_n$ maps the subalgebra $C^*_r(Z)\subset
C^*_u(S)$ into $\C[Z]\otimes M_{N_n}\subset l^\infty(X)\otimes M_{N_n}$ and
$\Psi_n$ maps $\C[Z]\otimes M_{N_n}$ into $C^*_r(Z)$. Then the nuclearity
of $C^*_r(Z)$ automatically follows.
So, let us examine the maps $\Phi_n,\Psi_n$. 
For each $n\geq 1$, we choose $N_n>0$ such that
$|B_{R_n}(S,x)|\leq N_n$ for all $x\in V(S)=X$.
Then, for each $x\in X$ we choose a subset $H^n_x\supset B_{R_n}(S,x)$ of size
  $N_n$
``locally''. That is, if $B_{S_n}(S,x)$ and $B_{S_n}(S,y)$ are rooted-labeled
isomorphic under the map $\theta$, then $\theta(H^n_x)=H^n_y$. Now for each $x\in
X$ let $P_n(x):l^2(X)\to l^2(H^n_x)$ be the orthogonal projection.
We set
$$\Phi_n:C^*_u(S)\to l^\infty(X)\otimes M_{N_n}(\C)$$
\noindent
by mapping $T$ to $\{P_n(x)T P_n(x)\}_{x\in X}$ in the same way as in
\cite{Roe}. The only difference between the approach of us and the one of
\cite{Roe} is the local choice of the projections $P_n$.
Clearly, if $T\in\C Z$ is a local kernel, then
$\Phi_n(T)\in\cA\otimes M_{N_n}(\C)$, where $\cA$ is the algebra defined in
Subsection \ref{s63}.
Hence, $\Phi_n$ maps the algebra $C^*_r(Z)$ into $\C[Z]\otimes M_{N_n}(\C)$.

\noindent
The maps $\Psi_n:l^\infty(X)\times M_{N_n}(\C)\to C^*_u(X)$ are defined by
mapping $\{T_x\}_{x\in X}$, $T_x\in B(l^2(H^n_x))\cong M_{N_n}(\C)$ to
$\sum_x M_n(x)^* T_x M_n(x)$, where $M_n(x)$ denotes the operator
of pointwise multiplication by the function $y\to \vs^n_y(x)$.
By the definition of Local Property A, the vectors $\vs^n_y$ are a priori
locally defined,
hence $\Psi_n$ maps $\cA\otimes M_{N_n}(\C)$ into $\C Z$. That is, $\Psi_n$
  maps
$\C[Z]\otimes M_{N_n}(\C)$ into $C^*_r(Z)$. Now our theorem follows. \qed
\subsection{Two examples for Local Property $A$}
\noindent
{\bf A tracial example.}  Let $Z$ be the generic URS constructed at the end of
Subsection \ref{coame1}.
That is, if $H\in Z$, then $S=S^Q_{\G_k}(H)$ is a colored graph satisfying
\begin{equation} \label{bi1}
C^1_\alpha r^\alpha\leq B_r(S,x)\leq C^2_\alpha r^\alpha\,
\end{equation}
\noindent
uniformly for some positive constants $C^1_\alpha$ and $C^2_\alpha$.
\begin{propo}
The graph $S$ has Local Property $A$.
\end{propo}
\proof
For a fixed vertex $w\in V(S)$ 
the unit vector $\vs^k_w$ is defined the following way.
\begin{itemize}
\item $\vs^k_w(z)=\frac{1}{\sqrt{|B_k(S,w)|}}$ if $z\in B_k(S,w)$.
\item $\vs^k_w(z)=0$ otherwise.
\end{itemize}
\begin{lemma} Let $x,y\in V(S)$ be arbitrary adjacent vertices. Then
$\|\vs^k_y-\vs^k_x\|^2\leq 2d\rho+\left(\frac{1}{\sqrt{2\rho d+1}} -1\right)^2\,.$
\end{lemma}
\proof Let $d$ be a bound for the vertex degrees of $S$ and let
$$\rho=\frac{|\partial B_k(S,x)|}{|B_k(S,x)|}\,.$$
We can suppose that $|B_k(S,x)|\leq |B_k(S,y)|$. Then we have that
$$\|\vs^k_y-\vs^k_x\|^2\leq 2d 
\frac{|\partial B_k(S,x)|}{|B_k(S,x)|}+|B_k(S,x)|\left(\frac{1}{\sqrt{|B_k(S,y)}}-
\frac{1}{\sqrt{|B_k(S,x)}}\right)^2\,.$$
\noindent
Now,
$|B_k(S,x)|\leq |B_k(S,y)|\leq (2\rho d+1)|B_k(S,x)|$.
Therefore,
$$\left(\frac{1}{\sqrt{|B_k(S,y)}}-
\frac{1}{\sqrt{|B_k(S,x)}}\right)^2\leq \frac{1}{|B_k(S,x)|}
\left(\frac{1}{\sqrt{2\rho d+1}} -1\right)^2\,,$$
\noindent
hence our lemma follows. \qed
\vskip 0.2in
\noindent
Thus, in order to prove our proposition, it is enough to show that for
every $\e>0$, there exists $K>0$ such that for each $x\in V(S)$
\begin{equation}\label{bi2}
\frac{|\partial B_k(S,x)|}{|B_k(S,x)|}\leq \e\,.
\end{equation}
\noindent
Note that (\ref{bi1}) implies that $S$ has the doubling condition, hence
(\ref{bi2}) follows from Theorem 4 \cite{Tess}. \qed
\vskip 0.2in
\noindent
{\bf A non-tracial example.}
Let $T$ be a $3$-regular tree. It is well-known that $T$ has Property $A$. The
construction goes as follows.
First, we pick an infinite ray $R=(x_0,x_1,\dots)$ towards the infinity. Then
for
each $t\in T$, there is a unique adjacent vertex $\phi(t)$ towards $R$ (if
$t=x_i$, $\phi(t)=x_{i+1}$). Then for a vertex $s$, we choose the path
$$P^n_s=(s,\phi(s),\phi^2(s),\dots, \phi^{n^2-1}(s))\,.$$
The unit vector $\vs^n_s$ is associated to the path $P^n_s$ as above, that is
$\vs^n_s(z)=\frac{1}{n}$ if $z\in P_s$ and $\vs^n_s(z)=0$ otherwise.
\begin{propo}
One can properly color $T$ by finitely many colors to obtain a Schreier graph
of Local Property $A$ that generates a generic URS.
\end{propo}
\proof
Our goal is to choose a coloring that encodes $\phi$. First, pick any
finite proper coloring $c:E(T)\to K$ for some finite set $K$
such that $c(e)\neq c(f)$ if $e\neq f$ and the distance
of $e$ and $f$ is less than $3$. Now we recolor the edge $(a,\phi(a))$
by $c(a,\phi(a))\times c(\phi(a),\phi^2(a))$. Hence, we obtained a proper
coloring $c':E(T)\to K\times K$ such that $\phi$ is encoded in the coloring so
the paths $P^n_s$ (and thus the unit vectors $\vs^n_s$) can be chosen locally.
Now let $m:E(T)\to A$ be the coloring given in Proposition \ref{constru}.
Then $m\times c':E(T)\to A\times K\times K$ provides a proper coloring of $T$,
such that the resulting Schreier graph has Local Property $A$ and generates a
generic $URS$. \qed

\section{The Feldman-Moore construction revisited}\label{feld}
Let $\alpha:\Gamma\to (X,\mu)$ be a  measure preserving action of a finitely generated group $\G$ 
on a standard probability measure space $(X,\mu)$.
The following construction is due to Feldman and Moore \cite{FM}.
We call a bounded measurable function $K:X\times X\to \C$ an $FM$-kernel if
\begin{itemize}
\item
$K(x,y)\neq 0$ implies that $x$ and $y$ are on the same orbit.
\item 
There exists a constant $w_K$ such that if $x$ and $y$ are on the orbit graph $S$, then
$d_S(x,y)> w_K$ implies that $K(x,y)=0\,.$
\end{itemize}
\noindent
The $FM$-kernels for the unital $*$-algebra $FM(\alpha)$, where
\begin{itemize}
\item $(K+L)(x,y)=K(x,y)+L(x,y)\,.$
\item $KL(x,y)=\sum_{z\in X} K(x,z) L(z,y)$\,.
\item $K^*(x,y)=\overline{K(y,x)}\,.$
\end{itemize}
The trace function $\Tr_\alpha$ is defined on $FM(\alpha)$ by
$$\Tr_\alpha(K)=\int_X K(x,x) d\mu(x)\,.$$
\noindent
Then, by the GNS-construction we can obtain a tracial von Neumann-algebra $FM(\alpha)\subset M(\alpha)$ in such a way that
the trace on $M(\alpha)$ is the extension of $\Tr_\alpha$. Let us very briefly
recall the construction. We define a pre-Hilbert space structure on
$FM(\alpha)$
by $\langle A,B\rangle=\Tr_{\alpha}(B^*A)$. Then $L_A(B)=AB$ defines
a map of $LM(\alpha)$ into $B(FM(\alpha))$. Then $M(\alpha)$ is
the weak closure of the image.  In particular, 
$\{K_n\}^\infty_{n=1}\subset M(\alpha)$ converges to $K\in M(\alpha)$ 
weakly if and only if for any $A,B\in LM(\alpha)$, $\lim_{n\to\infty}\Tr(AK_nB)=\Tr(AKB).$
Now, let $Z\subset\Sub(\Gamma)$ be a URS and $\mu$ be a $\G$-invariant Borel probability measure on $Z$. Again, $\beta$ denotes the $\G$-action on $Z$.
By definition, we have a natural homomorphism: $\phi_\beta:\C Z\to M(\beta)$.
\begin{propo}
The map $\phi_\beta$ is injective and $\phi_\beta(\C Z)$ is weakly dense in the von Neumann algebra $M(\beta)$.
Furthermore, the map $\phi_\beta$ extends to a continuous embedding $\overline{\phi}_\beta:C^*_r(Z)\to M(\beta)$.
\end{propo}
\proof First note, that if $K:X\times X\to \C$ is an $FM$-kernel, then $K$ can be written as $\sum^t_{i=1}M_{f_i} K_{g_i}$, where
\begin{itemize}
\item For any $1\leq i \leq t$, $f_i$ is a bounded $\mu$-measurable function.
\item $M_{f_i}\in FM(\beta)$ is supported on the diagonal and $M_{f_i}(x,x)=f_i(x)$.
\item $g_i\in\Gamma$ and $K_{g_i}(x,y)=1$ if $\beta(g_i)(y)=x$, otherwise $K_{g_i}(x,y)=0$.
\end{itemize}
\noindent
Let $0\neq K\in\C Z$. In order to prove that $\phi_\beta$ is injective, it is enough to show that
$\Tr_\beta(\phi_\beta(K^*K))\neq 0$.
Let 
$$U=\{x\in X\,\mid\, K^*K(x,x)\neq 0\}\,.$$
\noindent
Then $U$ is a nonempty open set, so by minimality of the action $\beta$, $\mu(U)>0$ since $\mu$ is $\Gamma$-invariant.
Therefore, $\Tr_\beta(\phi_\beta(K^*K))\neq 0$. Now we show that $\phi_\beta(\C Z)$ is weakly dense in the von Neumann algebra $M(\beta)$.
\begin{lemma} \label{weakconv}
Let $K_n\in FM(\beta)$, $K\in FM(\beta)$ such that
\begin{itemize}
\item
$\sup_{x,y\in X} |K_n(x,y)|<\infty\,,$
\item
$\sup_{n\geq 1} w_{k_n}<\infty\,,$
\item
For $\mu$-almost every $x$, $K_n(x,y)\to K(x,y)$ for all $y\in\Orb(x)$.
\end{itemize}
\noindent
then $\Tr_\beta(AKB)=\lim_{n\to\infty} \Tr_\beta(AK_n B)$ holds for any pair $A,B\in FM(\beta)$, hence by the GNS-construction
$\{K_n\}^\infty_{n=1}$ weakly converges to $K$.
\end{lemma}
\proof
Recall that
$$\Tr_\beta(AKB)= \int_X \sum_{y,z\in \Orb(x)} A(x,y) K(y,z) B(z,x) \,d\mu(x)\,.$$
\noindent
By our condition,
for almost every $x\in X$,
$$\lim_{n\to\infty} \sum_{y,z\in \Orb(x)} A(x,y) K_n(y,z) B(z,x)= \sum_{y,z\in \Orb(x)} A(x,y) K(y,z) B(z,x) \,,$$
hence by Lebesgue's Theorem
$$\Tr_\beta(AKB)=\lim_{n\to\infty} \Tr_\beta(AK_n B)\,.\quad \qed$$
\vskip 0.2in
\noindent
Now, let $K\in FM(\beta)$. We need to find a sequence $\{K_n\}^\infty_{n=1}\subset \phi_\beta(\C Z)$ that
weakly converges to $K$. Let $K=\sum^t_{i=1}M_{f_i}K_{\gamma_i}$. By Lemma \ref{alapok}, for
every $\alpha \in E_r$ we have a clopen set
$W_\alpha\subset Z$ such that
$\chi_{W_\alpha}\in \C Z$ and $\cup_{\alpha\in E_r} W_\alpha$ forms a partition of $Z$. Furthermore,
if $U\in Z$ is an open set, then we have a sequence
$\{Q^A_r\subset E_r\}$ so that
$$ \cup_{\alpha\in Q^A_1} W_\alpha \subset 
\cup_{\alpha\in Q^A_2} W_\alpha \subset \dots$$
\noindent
and
\begin{equation} \label{sec61}
\bigcup^\infty_{r=1}(\cup_{\alpha\in Q^A_r W_\alpha})=U\,.
\end{equation}
\noindent
Since $Z$ is homeomorphic to the Cantor set and $\mu$ is a Borel measure, for any $\mu$-measurable set $A\subset Z$ we have
a sequence of open sets $\{U_n\}^\infty_{n=1}\subset Z$ such that
\begin{equation} \label{sec62}
\{\chi_{U_n}\}^\infty_{n=1} \to \mu_A
\end{equation}
 $\mu$-almost everywhere.
Therefore, by (\ref{sec61}) and (\ref{sec62}), for any $1\leq i \leq t$,
we have a uniformly bounded sequence of functions
$\{g_{ij}\}^\infty_{j=1}$ tending to $f_i$ almost everywhere, such that for any
$i,j\geq 1$, $g_{ij}\in\cA$.
For $r\geq 1$, let 
$K_r=\sum^t_{i=1} M_{g_{ir}}K_{\gamma_i}\in \phi_\beta(\C Z)\,.$
Then for $\mu$-almost every $x\in X$
$$\lim_{r\to \infty} K_r(x,y)=K(x,y)$$
\noindent
provided that $y\in \Orb(x)$.
Therefore by Lemma \ref{weakconv}, $\{K_r\}^\infty_{r=1}$ weakly converges to $K$.
Hence, $\phi_\beta(\C Z)$ is weakly dense in $LM(\beta)$ and thus $\phi_\beta(\C Z)$ is weakly dense in $M(\beta)$ as well. 
Now we prove that $\phi_\beta$ extends to $C^*_r(Z)\,.$
First note that $\overline{Tr}_\beta(K)=\int K(x,x)\,d\mu(x)$ is a continuous trace on $C^*_r(Z)$ extending  $\Tr_\beta$.
Indeed, $\overline{Tr}_\beta(K)\leq \sup_{x\in X} |K(x,x)|\leq \|K\|\,.$
Let $\cN$ be the von Neumann algebra obtained from $C^*_r(Z)$ by the GNS-construction using the continuous trace $\overline{Tr}_\beta$.
The weak closure of $\C Z$ in $\cN$ is isomorphic to $M_\beta$, hence it is enough to prove that $\phi_\beta(\C Z)$ is weakly dense in
$C^*_r(Z)$. Let $A,B,K \in C^*_r(Z)$, $\{K_n\}^\infty_{n=1}\subset \C Z$, such that $K_n\to K$ in norm.
Then by the continuity of the trace, $\lim_{n\to\infty} \overline{Tr}_\beta(AK_n B)=  \overline{Tr}_\beta(AK B).$
Hence $\C Z$ is in fact weakly dense in $C^*_r(Z)$. \qed
\section{Coamenability and amenable traces} \label{amena}
\subsection{Amenable trace revisited}
First, let us recall the notion of amenable traces from \cite{Brown}.
Let $\cA$ be a $C^*$-algebra of bounded operators on the standard separable
Hilbert space $\cH$.
Let $\{P_n\}^\infty_{n=1}$ be a sequence of finite dimensional projections in
$\cH$ such that
\begin{itemize}
\item For any $a\in\cA$
$$\lim_{n\to\infty} \frac{\|AP_n-P_nA\|_{HS}}{\|P_n\|_{HS}}=0\,.$$
\item $$\tau(A)= \lim_{n\to\infty} \frac{\langle
  AP_n,P_n\rangle_{HS}}{\|P_n\|^2_{HS}}$$
\noindent
defines a continuous trace on $\cA$, where $\langle A,B
\rangle_{HS}=\Tr(B^*A)$ for Hilbert-Schmidt operators.
\end{itemize}
\vskip 0.2in
\noindent
Then $\tau$ is called an {\bf amenable trace}.
Now let $\Gamma$ be a finitely generated group as above and $Z\subset\Sub(\G)$ be a
coamenable generic URS. Let $H\in Z$, and consider the usual representation
of $C^*_r(Z)$ on $l^2(\Gamma/H)$ by kernels. Let $\{T_n\}^\infty_{n=1}$ be a
sequence
of induced subgraphs in $S=S^Q_\G(H)$ such that 
$\lim_{n\to\infty} \frac{|\partial T_n|}{|V(T_n)|}=0\,.$
Also, let us suppose that the sequence  $\{T_n\}^\infty_{n=1}$ is convergent
in the sense of Benjamini and Schramm as defined in Subsection \ref{sub43}.
Observe that convergence means that for any $r\geq 1$ and $\alpha\in E_r$
$$\lim_{n\to\infty}\frac{|V(T_n)\cap \alpha|}{|V(T_n)|}=t(\alpha)\,$$
\noindent
exists and $t(\alpha)=\mu(\lambda_H(\alpha))\,$ (see Subsection \ref{s63})\,,  where the $\G$-invariant
probability measure $\mu$ on $Z$ is
the limit of the sequence $\{T_n\}^\infty_{n=1}$. We define the amenable trace
$\tau$ similarly as in \cite{Elektrans}. For $n\geq 1$, let $P_n:l^2(\G/H)\to
l^2(V(T_n))\subset l^2(\G/H)$ be the orthogonal projection.
\begin{propo}
For any $K\in C^*_r(Z)$
$$\tau(K)=\lim_{n\to\infty} \frac{\langle
  AP_n,P_n\rangle_{HS}}{\|P_n\|^2_{HS}}$$
\noindent
exists and $\tau(K)=\overline{\Tr}_\mu(K)$ (as defined in Subsection
\ref{feld}).
Also, for any $K\in C^*_r(Z)$, 
$$\lim_{n\to\infty} \frac{\|KP_n-P_nK\|_{HS}}{\|P_n\|_{HS}}=0\,,$$
\noindent
hence $\tau$ is an amenable trace.
\end{propo}
\proof
Let us start with a simple observation.
\begin{lemma}\label{technical}
Let $\{H_n:\G/H\times \G/H\to\C\}^\infty_{n=1}$ be a sequence
of maps such that
\begin{itemize}
\item There exists $K>0$, $|H_n(x,y)|\leq K$, for any $n\geq 1$ and $x,y\in
  \G/H$.
\item $\lim_{n\to\infty}\frac{|Q_n|}{|V(T_n)|}=0\,,$ where
$$Q_n=\{(x,y)\in\G/H\times \G/H\,\mid\, H_n(x,y)\neq 0\}\,.$$
\end{itemize}
\noindent
Then $\lim_{n\to\infty} \frac{|\Tr(H_n)|}{\|P_n\|^2_{HS}}=0.$
\end{lemma}
\begin{lemma}
Let $K\in\C Z$. Then
$$\lim_{n\to\infty} \frac{\|KP_n-P_nK\|^2_{HS}}{\dim P_n}=0\,.$$
\end{lemma}
\proof
First we have that
$$\|KP_n-P_nK\|^2_{HS}=\Tr((P_nK^*-K^*P_n)(KP_n-P_nK))=\Tr(P_nK^*KP_n)-$$
$$- \Tr(K^*P_nKP_n)-\Tr(P_nK^*P_nK)+\Tr(K^*P_n K)\,.$$
\noindent
For $n\geq 1$, let $K_n:\G/H\to \G/H\to\C$ be defined in the following way.
$K_n(x,y)=K(x,y)$ if $x,y\in V(T_n)$, otherwise, $K_n(x,y)=0\,.$ That is, for
any $n\geq 1$, $K_n$ is a trace-class operator.
Now, we have that
$$\Tr(P_nK^*K P_n)=\Tr(P_nK_n^*K_nP_n)+\Tr(P_n(K-K_n)^* K_nP_n)+$$
$$+\Tr(P_nK^*(K-K_n)P_n)\,.$$
\begin{sublemma}
Both the sequences $$\{P_nK^*(K-K_n)P_n\}^\infty_{n=1} \quad \mbox{and} \quad
\{P_n(K-K_n)^*K_nP_n\}^\infty_{n=1}$$ \noindent  satisfy the conditions of our Lemma \ref{technical}.
\end{sublemma}
\proof
Notice that
$$(P_nK^*(K-K_n)P_n)(x,y)=\sum_{z\in \G/H} P_n(x,x)K^*(x,z)(K-K_n)(z,y)
P_n(y,y)\,.$$
\noindent
Observe that $(K-K_n)(z,y)P_n(y,y)\neq 0$ implies that $y\in V(T_n)$, $z\notin
V(T_n)$,
$d_{S^Q_\G(H)}(z,y)\leq w_K$. Since $\lim_{n\to\infty}
\frac{|\partial T_n|}{|V(T_n)|}=0\,,$ our sublemma immediately follows. \qed
\vskip 0.2in
\noindent
Repeating the arguments of our sublemma it follows that
$$\lim_{n\to\infty} \frac{\|KP_n-P_nK\|^2_{HS}}{\dim P_n}=$$
$$=\frac{\Tr(P_nK_n^*K_nP_n)-\Tr(K_n^*P_nK_nP_n)}{\dim P_n} +$$ $$+\frac{\Tr(K_n^*P_nP_nK_n)-\Tr(P_nK_n^*K_nP_n)}{\dim
  P_n}=0\,$$
\noindent
since $K_n$, $P_n$ are trace-class operators. \qed
\vskip 0.2in
\noindent
Now let $K\in C^*(Z)$ and $\e>0$.
Let $L\in \C Z$ such that $\|K-L\|<\e$. By the previous lemma, 
there exists $N>0$ such that if $n\geq N$, then $\|LP_n-P_nL\|_{HS}<\e\|P_n\|_{HS}$.
We have that
$$\|KP_n-LP_n\|_{HS}\leq\|K-L\| \|P_n\|_{HS}$$
\noindent
and
$$\|P_nK-P_nL\|_{HS}\leq\|K-L\| \|P_n\|_{HS}\,.$$
\noindent
Therefore, if $n\geq N$,
$\|KP_n-P_nK\|_{HS}<3\e\|P_n\|_{HS}$. Hence our proposition follows. \qed
\subsection{Uniformly amenable traces}
Recall from \cite{Brown} that an amenable trace $\tau$ is {\bf uniformly amenable}
if the resulting von Neumann algebra
$\pi_\tau(\cA)''$ is hyperfinite. Let $\alpha:\G\curvearrowright (Z,\mu)$ be
the Benjamini-Schramm limit of the graph sequence
$\{V(T_n)\}^\infty_{n=1}$ as in the previous subsection. Since
$C^*_r(Z)$ is dense in $M(\alpha)$, $\tau$ is a uniformly amenable trace
if and only if the equivalence relation generated by the action $\alpha$ is
hyperfinite.
By Theorem 1. \cite{Elekfun}, the equivalence relation above is hyperfinite if
and only if $\{T_n\}^\infty_{n=1}$ is a {\bf hyperfinite graph sequence}. That
is,
for any $\e>0$ there exists $K>0$ such that
for any $n\geq 1$ one can remove from the graph $T_n$ $\e|V(T_n)|$ edges 
in such a way that all the components of the remaining graph has at most $K$
elements.
Therefore we have the following proposition.
\begin{propo}\label{reflex}
Let $Z\subset\Sub(\Gamma)$ be a generic URS and $H\in Z$. If
$S^Q_\G(H)$ admits a convergent hyperfinite sequence of finite subgraphs
$\{T_n\}^\infty_{n=1}$
such that $\frac{|\partial T_n|}{|V(T_n)|}\to 0$, then $C^*_r(Z)$ has a
uniformly
amenable trace.
If
$S^Q_\G(H)$ admits a convergent nonhyperfinite sequence of finite subgraphs
$\{T_n\}^\infty_{n=1}$
such that $\frac{|\partial T_n|}{|V(T_n)|}\to 0$, then $C^*_r(Z)$ has a
non-uniformly
amenable trace, that is, by Theorem 4.3.3. $C^*_r(Z)$ is not locally
reflexive, hence it is a nonexact $C^*$-algebra.
\end{propo}

\section{A nonexact example}
\subsection{The construction}
In Section 4.3.3 of \cite{Cecc}, we constructed a Schreier graph
$S=S^Q_{\G}(H)$ of a group $\G$ such that
that the orbit closure of $S$ is a generic URS. Also, $S$ contains
a convergent, hyperfinite sequence of finite subgraphs $\{T_n\}^\infty_{n=1}$
with $\frac{|\partial T_n|}{|V(T_n)|}\to 0$ and a convergent nonhyperfinite
sequence of finite subgraphs $\{W_n\}^\infty_{n=1}$
with $\frac{|\partial W_n|}{|V(W_n)|}\to 0$. Let $Z$ be the orbit closure of
$H$ in $\Sub(\G)$. Then by Proposition \ref{reflex}, $C^*_r(Z)$ is not a
locally reflexive (hence nonexact) simple, unital separable $C^*$-algebra
with both uniform amenable and non-uniform amenable traces.
For completeness, we present a somewhat slicker construction, that is very
similar to the one given in \cite{Cecc} 
\vskip 0.1in
\noindent
{\bf Step 0.} For $n\geq 1$, let $C_n$ be the cycle of length $2^{n+1}$. Let
$\G_3$ be the free group of three cycle groups of rank $2$. Let $\G_3\subset
M_1 \supset M_2 \supset \dots, \cap^\infty_{n=1} M_n=\{1\}$ be a sequence of
finite index normal subgroups. For $n\geq 1$ let $V_n$ be the
underlying graph of the Cayley-graph of $\G_3/M_n$. That is, the sequence
$\{V_n\}^\infty_{n=1}$ itself is a convergent nonhyperfinite sequence of
finite graphs. 
\vskip 0.1in
\noindent
{\bf Step 1.} Let $G_1=C_1$ and $H_1=C_2$.
\vskip 0.1in
\noindent
{\bf Step n.} Let us suppose that the graphs $G_2\subset G_4\subset \dots\subset
G_{2n-2}$, $G_1\subset
G_3\subset\dots\subset G_{2n-1}$ 
and $\{H_i\}^{2n-1}_{i=1}$ are already defined in such a way that
\begin{itemize}
\item If $1\leq i \leq n-1$, then $H_{2i}=V_{k_i}$ for some $k_i>0$,
$H_{2i+1}=C_{l_i}$ for some $l_i>0$.
\item
For any $1\leq i \leq n-1$
we have disjoint subsets \\
$\{R^j_{2i}\}^i_{j=1}\subset H_{2i}, \{R^j_{2i+1}\}^i_{j=1}\subset H_{2i+1}$
so that
$$|R^j_{2i}||V(G_j)|\leq \frac{1}{10^{j}} |V(H_{2i})|\,\,\,,\,\,\,
|R^j_{2i+1}||V(G_j)|\leq \frac{1}{10^{j}} |V(H_{2i+1})|\,.$$
\item
We have positive integers $T_1< T_2<\dots T_{n-1}$ such that
for any $1\leq j\leq i\leq n-1$
$$ R^j_{2i}\cap B_{T_j}(H_{2i},x)\neq \emptyset, R^j_{2i+1}\cap
  B_{T_j}(H_{2i+1},y)\neq\emptyset\,,$$
for all $x\in V(H_{2i}), y\in V(H_{2i+1})\,.$
\item The graph $G_{2i}$ is constructed in such a way that for any 
$1\leq j\leq i$  
a copy of $G_j$ is connected to all the vertices of $R^j_{2i}\subset
V(H_{2i})$.
The graph $G_{2i+1}$ is constructed in such a way that for any 
$1\leq j\leq i$  
a copy of $G_j$ is connected to all the vertices of $R^j_{2i+1}\subset
V(H_{2i+1})$. Connecting a graph $A$ to a graph $B$ means that we add a
disjoint copy of $A$ to $B$ plus an extra edge beween a vertex of $A$ and a
vertex of $B$.
\end{itemize}
\noindent
Now we construct the graphs $G_{2n}$ and $G_{2n+1}$. We pick a graph
$H_{2n}=V_{k_n}$ and $H_{2n+1}=C_{l_n}$ in such a way that
$$|V(G_n)|\leq \frac{1}{10^n} |V(H_{2n})|\,\,,\,\,
|V(G_n)|\leq \frac{1}{10^n} |V(H_{2n+1})|\,.$$
\noindent
We define $T_n$ as the maximum of the diameters of the graphs $H_{2n}$ and
$H_{2n+1}$.
Now we use the fact that we have normal covering maps $\zeta_{2n}:H_{2n}\to
H_{2n-2}$
and $\zeta_{2n+1}:H_{2n+1}\to H_{2n-1}$.
For $1\leq j\leq n-1$, set $R^j_{2n}=\zeta^{-1}_{2n}(R^j_{2n-2})$ and
$R^j_{2n+1}=\zeta^{-1}_{2n+1}(R^j_{2n-1})$.
That is, 
$$|R^j_{2n}||V(G_j)|\leq \frac{1}{10^j} |V(H_{2n})|\,\,\,,\,\,\, 
|R^j_{2n+1}||V(G_j)|\leq \frac{1}{10^j} |V(H_{2n+1})|\,.$$
Now for any $1\leq j \leq n-1$, we connect a copy of $G_j$ to the vertices of
$R^j_{2n}$
and a copy of $G_j$ to the vertices of
$R^j_{2n+1}$. Finally, we connect a copy of $G_n$ to a single vertex $a$ of $V(H_{2n})$
not covered by any $R^j_{2n}$ and a copy of $G_n$ to a single vertex $b$ of $V(H_{2n+1})$
not covered by any $R^j_{2n+1}$. Then we set $a=R^n_{2n}$, $b=R^n_{2n+1}$.
Then, for any $1\leq j \leq n$, $R^j_{2n}\cap B_{T_j}(H_{2n},x)\neq \emptyset$
and $R^j_{2n+1}\cap B_{T_j}(H_{2n+1},y)\neq \emptyset$, for all $x\in V(H_{2n})$
and $y\in V(H_{2n+1})$. Also, we have that 
$$|R^n_{2n}||V(G_n)|\leq \frac{1}{10^n} |V(H_{2n}|\,\,\,,\,\,\, 
|R^n_{2n+1}||V(G_n)|\leq \frac{1}{10^n} |V(H_{2n+1}|\,.$$
\noindent
Now, we have graphs $G_1\subset  G_3\subset G_5\subset\dots$ and set
$G=\cup_{i=1}^\infty G_{2i-1}$. 
By the self-similar nature of our construction, it is easy to see that for any
$x\in V(G)$ and $r\geq 1$, there exists $N_{x,r}>0$ such that
if $y\in V(G)$, then there is a $z\in B_{N_{x,r}}(G,y)$ such that
$B_r(G,x)$ and $B_r(G,z)$ are isomorphic as rooted graphs.
Notice that $G_1\subset G_3\subset\dots$ are forming a hyperfinite sequence of
subgraphs.
so that $\frac{|\partial G_{k_i}|}{|V(G_{k_i})|}\to 0$.
Also, we have subgraphs $G_2$, $G_4$,\dots connected by one single edge to
$H_{2n+1}$, that are forming such a sequence of finite graphs that no
subsequence of them is hyperfinite and $\frac{|\partial G_{2i}|}{|V(G_{2i})|}\to 0$. Indeed, in the construction $H_{2n}\subset
G_{2n}$, $|V(H_{2n})|> |V(G_{2n})|/2$ and the sequence of finite graphs
$\{H_{2n}\}^\infty_{n=1}$ is a large girth graph sequence so no subsequence of
$\{G_{2n}\}^\infty_{n=1}$ 
can be hyperfinite (since that would mean that a subsequence of
$\{H_{2n}\}^\infty_{n=1}$ is hyperfinite as well). 
Now we can apply the coloring construction of Proposition \ref{constru} to
obtain the colored graph $S$ (and hence a URS $Z\subset\Sub(\G_k)$ we are sought
of. In the graph $S$ there is a convergent hyperfinite sequence of finite
subgraphs $\{T_n\}^\infty_{n=1}$ with $\frac{|\partial T_n|}{|V(T_n)|}\to 0$ and
there is a convergent nonhyperfinite sequence of finite
subgraphs $\{T'_n\}^\infty_{n=1}$ with $\frac{|\partial T'_n|}{|V(T'_n)|}\to 0$

\subsection{Two more interesting properties of the nonexact URS}
Let $Z\subset\Sub(\G)$ be the generic URS constructed above.
\begin{propo}
There is no free continuous action of any countable group that is 
Borel orbit equivalent
to $Z$.
\end{propo}
\proof Suppose that $Z$ is Borel orbit equivalent to a free continuous
action $\theta$ of the countable group $\Delta$. If $\Delta$ is amenable, then 
any invariant measure of the action $\theta$ makes the equivalence relation
hyperfinite. However, $Z$ admits invariant measure that makes it a
nonhyperfinite measurable equivalence relation. If $\Delta$ is nonamenable
then all the invariant measures make the equivalence relation of 
the action $\theta$ nonhyperfinite. However, $Z$ has an invariant measure that
makes its equivalence relation hyperfinite (see also \cite{HM} for
a minimal action of a group that is not Borel equivalent to a free 
continuous action). \qed
\vskip 0.2in
\noindent
If $K_1, K_2\in Z$ then their Schreier graphs are locally indistinguishable. 
However,
it is possible that their Schreier graphs look globally quite different.
\begin{propo}
There exists $K_1,K_2\in Z$ such that $S^Q_\G(K_1)$ is one-ended and
$S^Q_\G(K_2)$ is two-ended.
\end{propo}
\proof
Clearly, there must be an element $K_1$ in $Z$ such that
the underlying graph of $S^Q_\G(K_1)$ is isomorphic to the graph $G$ in
our construction and $G$ is clearly one-ended. Now we pick a sequence of points
$x_n\in H_{2n+1}$.  Consider
a limitpoint $T=S^Q_\G(K_2)$ of the rooted Schreier graphs
$\{(S^Q_\G(K_1),x_n)\}^\infty_{n=1}$ in the compact space $\Sch^Q_\G$.
It is easy to see that $T$ has multiple ends. \qed

\vi
\vi
G\'abor Elek

\noindent
Lancaster University

\noindent
g.elek@lancs.ac.uk
\end{document}